\documentclass[11pt,a4paper,twoside]{amsart}
\usepackage{amsmath,amssymb,amsfonts,amsthm}
\usepackage{hyperref,color}
\usepackage{times}
\usepackage{graphicx}

\newtheorem{The}{Theorem}[section]

\newtheorem{Proposition}[The]{Proposition}
\newtheorem{Lemma}[The]{Lemma}
\newtheorem{Corollary}[The]{Corollary}

\theoremstyle{definition}
\newtheorem{Definition}[The]{Definition}
\newtheorem{Remark}[The]{Remark}
\newtheorem{Example}[The]{Example}


\begin{document}

\title{Lie algebras of infinitesimal CR-automorphisms
\\
of finite type, holomorphically nondegenerate,
\\
weighted homogeneous CR-generic submanifolds of $\mathbb{ C}^N$}

\author{Masoud Sabzevari}
\address{Department of Pure Mathematics,
University of Shahrekord, 88186-34141 Shahrekord, IRAN and School of
Mathematics, Institute for Research in Fundamental Sciences (IPM), P.O.Box: 19395-5746, Tehran, IRAN}
\email{sabzevari@math.iut.ac.ir}

\author{Amir Hashemi}
\address{Department of Mathematical Sciences,
Isfahan University of Technology, Isfahan, IRAN and School of Mathematics, Institute for Research in Fundamental
Sciences (IPM), Tehran, P.O.Box: 19395-5746, IRAN}
\email{amir.hashemi@cc.iut.ac.ir}

\author{Benyamin M.-Alizadeh}
\address{School of Mathematics and Computer Sciences, Damghan University, P.O.Box 3671641167, Damghan, IRAN}
\email{benyamin.m.alizadeh@gmail.com}

\author{Jo\"el Merker}
\address{D\'epartment de Math\'ematiques d'Orsay, B\^{a}timent 425,
Facult\'{e} des Sciences, Universit\'{e} Paris XI - Orsay, F-91405
Orsay Cedex,  FRANCE} \email{merker@dma.ens.fr}

\date{\number\year-\number\month-\number\day}

\maketitle

\begin{abstract}
We consider the significant class of holomorphically nondegenerate
CR manifolds of finite type that are
represented by some weighted homogeneous
polynomials and we derive some useful features which enable us
to set up a fast effective algorithm to compute their Lie
algebras of infinitesimal CR-automorphisms. This algorithm
mainly relies upon a natural
gradation of the sought Lie algebras, and it also
consists in treating {\em separately} the related graded
components. While some other methods are
based on constructing and solving an
associated {\sc pde} systems which
become time consuming as soon as the number of variables increases, the new
method presented here is based on plain techniques of linear algebra.
Furthermore, it benefits from a {\it divide-and-conquer} strategy to
break down the computations into some simpler sub-computations.
Moreover, we consider the new and effective concept of comprehensive Gr\"obner systems which provides us
some powerful tools to treat the computations in the parametric cases. The designed algorithm is also implemented in the {\sc Maple} software.
\end{abstract}

\pagestyle{headings} \markright{Lie algebras of infinitesimal
CR-automorphisms }

\section{Introduction}
\label{Introduction}

Let $M\subset\mathbb C^{n+k}$ be a Cauchy-Riemann (CR for short)
submanifold of CR dimension $n \geqslant 1$ and of codimension $k
\geqslant 1$ ({\it see} $\S$\ref{General description} for all
pertinent definitions used in this introduction) represented in
coordinates $z_j$ and $w_l:=u_l+iv_l$ for $j=1,\ldots,n$ and
$l=1,\ldots,k$. As is standard in CR-geometry (\cite{BER, Boggess,
Beloshapka2004, Beloshapka2002}), one may often assign {\it weight}
$[z_j] := 1$ to all the complex variables $z_j$ and weights some
certain integers $[w_l]\in\mathbb N$ with $1<[w_1]\leqslant
[w_2]\leqslant\cdots [w_k]$ to the variables $w_1, \dots, w_k$.
Accordingly, the weight of the conjugation of each complex variable
and of its real and imaginary parts as well are all equal to that of
the variable and, moreover, the assigned weight of any constant
number $a\in\mathbb C$ and of coordinate vector fields are:
\begin{equation*}
\aligned
&
[a]:=0, \ \ \ \ [\partial_{z_j}]:=-[z_j], \ \ \ \
[\partial_{w_l}]:=-[w_l], \ \ \ \ \ \ {\scriptstyle (j=1,\ldots,n, \
\ \ l=1,\ldots,k)}.
\endaligned
\end{equation*}
Furthermore, the weight of a monomial $F$ in the $z_j$, $w_k$,
$\overline{ z}_j$, $\overline{ w}_k$ is the sum
taken over the weights of all variables of $F$ with regards to their
powers and also, the weight of each coordinate vector field of the form
$F\,\partial_x$ with $x=z_i,w_l$ is defined as $[F]-[x]$. For
instance, we have $[(a\,z\,w_1^2)\,\partial_{w_2}]=[z]+2[w_1]-[w_2]$. A
polynomial or a vector field is called {\sl weighted homogeneous of
weight} a certain
integer $d$ whenever each of its terms is of homogeneity $d$.

As is well-known ({\it see} \cite{BER} Theorem 4.3.2 and the Remark
after it or \cite{Merker-Porten-2006} Theorem 2.12), every real analytic
generic CR manifold $M$ of CR dimension $n$ and codimension $k$ can be
represented locally in a neighborhood of the origin as the graph of
$k$ defining equations of the form:
\begin{equation}
\aligned \label{model-1}
\left\{
\begin{array}{l}
v_1:=\Phi_1(z,\overline z, u)+{\rm o}([w_1]),
\\
\ \ \ \ \ \ \ \vdots
\\
v_k:=\Phi_k(z,\overline z, u)+{\rm o}([w_k]),
\end{array}
\right.
\endaligned
\end{equation}
where the weight of all variables $z_j$ is $1$ and
the weights of the variables
$w_l$ are the H\"ormander numbers of $M$ at the origin. Moreover, each function
$\Phi_l$ is a weighted homogeneous polynomial of weight $[w_l]$ and
${\rm o}(t)$ denotes
remainder terms having weights $> t$.

For a CR manifold $M$ passing through the origin, the Lie group ${\sf
Aut}_{CR}(M)$ is the holomorphic symmetry group of $M$, that is the
local Lie
group of local biholomorphisms mapping $M$ to itself.
The Lie algebra $\frak{aut}_{CR}(M)$, associated to this group is
called the {\sl Lie algebra of infinitesimal CR-automorphisms} of $M$
and it consists of all holomorphic vector fields\,\,---\,\,$(1,0)$
fields with holomorphic coefficients\,\,---\,\,whose real parts are
tangent to $M$. Due to the fact that many geometric features of
CR manifolds can be investigated by means of their associated Lie
algebras of infinitesimal CR-automorphisms and because of
central applications in Cartan geometry and in Tanaka
theory exist, studying such algebras have gained an increasing interest
during the recent decades ({\it cf.} \cite{Beloshapka2004, BES,
5-cubic,MS, Stanton}). As is known, the Lie algebra
$\frak{aut}_{CR}(M)$ is finite dimensional if and only if $M$ is
holomorphically nondegenerate and of finite type (\cite{BER,Stanton,
Gaussier-Merker,
Juhlin,Merker-Porten-2006}).

Consider the complex space $\mathbb C^{n+k}$ equipped with the
coordinates $z_1,\ldots,z_n,w_1,\ldots,w_k$, where $w_j:=u_j+iv_j$,
assume again that certain weights have been assigned, and consider
holomorphically nondegenerate weighted homogeneous CR manifolds $M
\subset \mathbb C^{n+k}$ represented as graphs of $k$ certain
polynomials:
\begin{equation}
\aligned \label{model} M:=\left\{(z,w): \ \ \left[
\begin{array}{l}
\Xi_1(v_1,z,\overline z,u):=v_1-\Phi_1(z,\overline z,u)\equiv 0,
\\
\Xi_2(v_2,z,\overline z,u):=v_2-\Phi_2(z,\overline z,u)\equiv 0,
\\
\ \ \ \ \ \ \ \vdots
\\
\Xi_j(v_j,z,\overline z,u):=v_j-\Phi_j(z,\overline z,u)\equiv 0,
\\
\ \ \ \ \ \ \ \vdots
\\
\Xi_k(v_k,z,\overline z,u):=v_k-\Phi_k(z,\overline z,u)\equiv 0,
\end{array}
\right.\right\},
\endaligned
\end{equation}
with the right-hand sides $\Phi_j$ being weighted homogeneous
polynomials of weight equal to that $[v_j] = [w_j]$ of the left-hand
sides.  We will also assume that such $M$ are {\sl homogenous} in Lie
theory's sense, namely that ${\rm Aut}_{ CR} ( M)$ is locally
transitive near the origin.  Since it may arise some confusion with
the `homogeneous spaces' terminology, let us stress that we will
always use {\sl weighted homogeneous} about functions, and plainly
{\sl homogeneous} about ${\rm Aut}_{ CR} ( M)$.

The so introduced class of CR-generic manifolds is already an
extremely wide class in CR-geometry which includes of course well
known quadric CR-models such as those of Poincar\'{e} \cite{Poincare}
or of Chern-Moser \cite{Chern-Moser}, and also, there is nowadays an
extensive literature dealing with constructing a great number of such
weighted homogeneous CR manifolds ({\it see} for example
\cite{Baouendi, BER, Bloom, Mamai2009,Shananina2000} and
\cite{Beloshapka2012}--\cite{Beloshapka2001}), their Lie groups being
far from being completely understood.

Indeed, in a series of recent papers (for instance
\cite{Beloshapka2012}--\cite{Beloshapka1997}), Valerii Beloshapka
studied extensively the subject of model surfaces and found some
considerable results in this respect. Specifically in
\cite{Beloshapka2004}, he introduced and established the structure of
some nondegenerate models associated (uniquely) to totally
nondegenerate germs having arbitrary CR dimensions and real codimensions. He
also developed systematic tools for their construction. All of these models
are again included in our already mentioned class of CR manifolds.
Each of Beloshapka's model $M\subset\mathbb C^{n+k}$ of certain CR
dimension and codimension $n$ and $k$ enjoys some {\it nice}
properties (\cite{Beloshapka2004}, page 484, Theorem 14) which have
been encouraging enough to merit further investigation. In
particular, computing their Lie algebras of infinitesimal
CR-automorphisms and studying their structures may reveal some
interesting features of these CR-models, and also of all totally
nondegenerate CR manifolds corresponding to them ({\it cf.}
\cite{Beloshapka2011}).

It is worth noting that, traditionally the subject of computing Lie
algebras of infinitesimal CR-automorphisms is concerned with
expensive computations and the cost of calculations increases as much
as the number of the variables\,\,---\,\,namely the dimension or
codimension of the CR manifolds\,\,---\,\,increases. Indeed, solving
the {\sc pde} systems arising during these computations ({\it see}
subsection \ref{SHAM}) forms the most complicated part of the
procedure. That may be the reason
why, in contrast to the importance of the
subject, the number of the relevant computational works is still
limited (one finds some of them in \cite{MS, Beloshapka1997,
Shananina2000, Mamai2009}).

Very recently, the authors provided in \cite{HAMS} a new
general algorithm to compute the desired algebras by means of the
effective techniques of differential algebra ({\it see}
$\S$\ref{General description} below for a brief description of the
algorithm). It enables one to use conveniently the ability of
computer algebra for managing the associated computations of the
concerning {\sc pde} systems. Although this (general) algorithm is
able to decrease a lot the complexity of the computations and in
particular to utilize systematically the ability of computer algebra,
but {\it because of dealing with the {\sc pde}
systems}\,\,---\,\,which are complicated in their
spirit\,\,---\,\,the computations are still expensive in essence.

In the present paper
we aim to study, by means of a {\it weight analysis approach and with
an algorithmic treatment}, the intrinsic properties of the under
consideration CR manifolds in order to provide an effective
algorithm\,\,---\,\,entirely different and more powerful than that of
\cite{HAMS}\,\,---\,\,to compute the associated Lie algebras of
infinitesimal CR-automorphisms. The results enable us to bypass
constructing and solving the arising systems of {\sc pde} (which is
the classical method of
\cite{Beloshapka1997,MS,5-cubic,Mamai2009,Shananina2000, HAMS}) and
to reach the sought algebras by employing just simple techniques of
linear algebra.  This  decreases considerably the cost of computations,
hence simultaneously increases the performance of the algorithm.

We show that for a homogeneous CR manifold $M$ represented as
\thetag{\ref{model}}, the sought algebra $\frak
g:=\frak{aut}_{CR}(M)$ takes the graded form (in the sense of
Tanaka):
\begin{equation}
\label{g-int}
 \frak g=\underbrace{\frak g_{-\rho}\oplus\cdots\oplus\frak
g_{-1}}_{\frak g_-}\oplus\frak g_0\oplus\underbrace{\frak
g_1\oplus\cdots\oplus\frak g_\varrho}_{\frak g_+} \ \ \ \ \
\rho,\,\varrho\in\mathbb N.
\end{equation}
where each component $\frak g_t$ is the Lie subalgebra of all
weighted homogeneous vector fields having weight $t$. Summarizing, the
results provide us with the ways of:

\begin{itemize}

\smallskip\item[$\bullet$]
using a {\it divide-and-conquer} method to break down the
computations of the sought (graded) algebra $\mathfrak{ g}$ into some simpler
sub-computations of its (homogeneous) components $\frak g_t$;

\smallskip\item[$\bullet$]
employing some simple techniques of linear algebra
for computing these components $\mathfrak{ g}_t$
of infinitesimal CR-automorphisms
{\em without} relying on solving the {\sc pde} systems;

\end{itemize}\smallskip

We also find a simple criterion for how long it is necessary to
compute the homogeneous components $\mathfrak{ g}_t$ of $\frak g$. In
other words, it supplies us to recognize\,\,---\,\,in an algorithmic
point of view\,\,---\,\,where the {\it maximum homogeneity $\varrho$}
in \thetag{\ref{g-int}} is, while the minimum
homogeneity $\rho$ is easy to determine.
This criterion fails when ${\rm Aut}_{ CR}
(M)$, or equivalently $\mathfrak{ aut}_{ CR} ( M)$, is not locally
transitive. However without local transitivity, one still is able to
compute the homogeneous components $\frak g_t$, separately.

One of the main\,\,---\,\,somehow hidden\,\,---\,\,obstacles appearing among the computations
arises when the set of defining equations includes some certain {\it parametric} polynomials.
This case is quite usual as one observes in \cite{Beloshapka2004,Beloshapka1997,Mamai2009,Shananina2000}. To treat such cases, we suggest the modern and effective concept of {\sl comprehensive Gr\"obner systems} \cite{Weis, kapur,kapur2,monts1,monts3} which enables us to consider and solve (linear) parametric systems appearing among the computations.

This weight analysis computational approach to the determination of
infinitesimal CR-automorphisms can provide useful tools to study some
general problems in CR-geometry, as for instance the models of Kolar
(\cite{Kolar}), and also, it can offer numerical evidence to support a
recent conjecture of Beloshapka \cite{Beloshapka2012} which asserts
that:

\begin{itemize}
\item[]{\bf Beloshapka's Maximum conjecture:}
For a homogeneous CR-model as those of
\cite{Beloshapka2004} with the graded Lie algebra of infinitesimal
CR-automorphisms like \thetag{\ref{g-int}}, if $3\leqslant\rho< \infty$
then the corresponding universal CR manifold has `{\sl rigidity}' {\em
i.e.} $\varrho=0$.
\end{itemize}

This paper is organized as follows. Section \ref{General description}
presents a brief description of required very basic
definitions, notations and
terminology and also a presentation of our recently prepared algorithm
\cite{HAMS}. In Section \ref{section
rigidity-similarity}, we inspect some of the Lie algebras computed in
\cite{Beloshapka1997, MS, Shananina2000, Mamai2009} and we
observe a striking {\it similarity} in the coefficients of the
appearing holomorphic vector fields. This leads us to discover the
key entrance to the desired algorithm. In Section
\ref{section computation the components}, we employ the results of
the previous section to provide the strategy of computing separately
the homogeneous components of the sought algebras. We also provide
the necessary criterion for terminating such computations.
Section \ref{section main strategy} presents the desired algorithm of
computing the Lie algebras of infinitesimal CR-automorphisms
associated to the universal CR manifolds by using the results
obtained in the earlier sections. Finally, in section \ref{Comprehensive} we introduce briefly the modern concept of comprehensive Gr\"obner systems and show how it provides some effective tools to consider and solve appearing (linear) parametric systems.

The algorithm designed in this paper is implemented in {\sc Maple} 15 as the library {\sc CRAut}, accessible online as \cite{Implementation}. To do this, at first we needed to implement the algorithm {\sc PGB} introduced in the recently published paper \cite{kapur2} which enables us to consider the parametric defining equations in {\sc CRAut}.

\section{Basic preliminaries and definitions}
\label{General description}

On an arbitrary even-dimensional real vector space $V$, a {\sl
complex structure map} $J:V\rightarrow V$ is an $\mathbb{ R}$-linear map
satisfying $J\circ J=-{\rm Id}_V$. For example, in the simple case
$V:=T_p\mathbb R^{2N}=T_p\mathbb C^N, \ N\in\mathbb N$ with the local
coordinates $(z_1:=x_1+i\,y_1,\ldots,z_N:=x_N+i\,y_N)$ and with
$p\in\mathbb C^N$, one defines the (standard) complex structure by:
\[\footnotesize\aligned
J\Big(\frac{\partial}{\partial x_j}\Big)
=
\frac{\partial}{\partial
y_j}, \ \ \ \ \ \ J\Big(\frac{\partial}{\partial
y_j}\Big)=-\frac{\partial}{\partial x_j}, \ \ \ \  \ j=1,\ldots,N.
\endaligned
\]
One should notice that in general, an arbitrary subspace $H_p$ of
$T_p\mathbb C^N$ is not invariant under the complex structure map
$J$. Thus, one may give special designation to the largest
$J$-invariant subspace of $H_p$ as:
\[
H_p\cap J(H_p)=:H_p^c,
\]
which is called the {\sl complex tangent subspace of} $H_p$. Due to the
equality $J\circ J=-{\rm Id}$, this space is even-dimensional, too.

Similarly, one also introduces the smallest $J$-invariant real subspace
of $T_p\mathbb C^N$ which {\em contains} $H_p$:
\[
H^{i_c}_p:=H_p+J(H_p)
\]
and one calls its
the {\sl intrinsic
complexification} of $H_p$.

As an application,
consider the linear subspace $H_p:=T_pM$ of $\mathbb C^N$, for
some arbitrary connected differentiable submanifold $M$ of $\mathbb
C^N$. In general, it is not at all true that the complex-tangent
planes:
\[
T_p^cM=T_pM\cap J(T_pM)
\]
have constant dimensions as $p$ varies in $M$.

\begin{Definition}
Let $M$ be a real analytic submanifold of $\mathbb C^N$. Then $M$ is
called {\sl Cauchy-Riemann} (CR for short), if the complex dimensions
of $T^c_pM$ are {\em constant} as $p$ varies on $M$. Furthermore,
$M$ is
called {\sl generic} whenever:
\[
T_p^{i_c}M:=T_pM+J(T_pM)=T_p\mathbb C^N
\]
for each $p\in M$, which implies that $M$ is CR thanks
to elementary linear algebra.
\end{Definition}

One knows (\cite{ Boggess, BER, Merker-Porten-2006})
that any CR real analytic $M \subset \mathbb{ C}^N$
is contained in a thin strip-like complex submanifold $\mathcal{ C}'
\cong \mathbb{ C}^{ N'}$
with $N' \leqslant N$
in which $M \subset \mathbb{ C}^{ N'}$
is CR-generic, hence there is no restriction to
assume that $M$ is CR-generic, as we will always do here,
since we are interested in local CR geometry.

For such a CR-generic manifold $M$, we call the complex dimension of $T_p^cM$
by the {\sl CR-dimension} of $M$. Moreover, the subtraction of the
real dimension of $T_p^cM$ from that of $T_pM$, which is in fact the
real dimension of the so-called {\sl totally real part} $T_pM/T^c_pM$
of $T_pM$, coincides then with the
real {\sl codimension} of $M \subset \mathbb{ C}^N$.

A tangent vector field:
\[
{\sf X}:=\sum_{j=1}^N\,\bigg(a_j\,\frac{\partial}{\partial
Z_j}+b_j\,\frac{\partial}{\partial \overline Z_j}\bigg)
\]
of the complexified space $\mathbb CT_p\mathbb C^N:=\mathbb C\otimes
T_p\mathbb C^N$ is called of type
$(1,0)$ whenever all $b_j\equiv 0$
and is called of type $(0,1)$ whenever all $a_j\equiv 0$.
One denotes by $T^{1,0}_p\mathbb C^N$ and
$T^{0,1}_p\mathbb C^N$ the corresponding
subspaces of $\mathbb{ C}
T_p \mathbb{ C}^N$. Accordingly, for a CR manifold $M$ and an arbitrary
point $p$ of it, we denote:
\[
T^{1,0}_pM:=T^{1,0}_p\mathbb C^N\cap\mathbb C T_pM, \ \ \ \ {\rm and}
\ \ \ \  T^{1,0}_pM:=T^{0,1}_p\mathbb C^N\cap\mathbb C T_pM.
\]
One easily verifies the equality $T^{0,1}_pM=\overline{T^{1,0}_pM}$.
It is proved that ({\it see} \cite{Baouendi}, Proposition 1.2.8) the
complex tangent space $T^c_pM$ is the real part of $T^{1,0}_pM$, {\it
i.e.} $T^c_pM=\{{\sf X}+\overline{\sf X}: \, {\sf X}\in
T^{1,0}_pM\}$. Moreover, the complexified space $\mathbb C T^c_pM$ is
equal to the direct sum $T^{1,0}_pM\oplus T^{0,1}_pM$.

If $n$ and $k$ are the CR dimension and the real
codimension of a real analytic CR-generic submanifold $M\subset\mathbb
C^N$, then of course $N=n+k$, and also $M$ can be represented (locally)
by $k$ real analytic graphed equations:
\begin{equation}\label{Phi} {\rm
Im}\,w_j:=\Phi_j(z,\overline z,{\rm Re}\,w), \ \ \ \ \ \ \ \
{\scriptstyle (j=1,\ldots,k)}
\end{equation}
with some real-valued defining functions $\Phi_\bullet$ enjoying
the {\sl no-pluriharmonic term} condition:
\[
0\equiv\Phi_\bullet(z,0,{\rm Re}\,w)=\Phi_\bullet(0,\overline z,{\rm
Re}\,w).
\]
Solving the above real-valued defining functions in $w$ or in $\overline
w$, one also reformulates the defining functions of $M$ as
complex defining equations of the kind:
\begin{equation}
\label{Xi} w_j+\overline{w}_j=\overline{\Xi}_j(z,\overline{z},w), \ \
\ \ \ \ \ \ \ {\scriptstyle (j=1,\ldots,k)}.
\end{equation}
The real analyticity of the functions $\Xi_j$ enables us to employ their
Taylor series for computing the associated algebras of infinitesimal
CR-automorphisms ({\it see} the next subsection).

For a real analytic CR-generic manifold $M$ represented by the above
defining functions as \thetag{\ref{Xi}} and for each $p\in M$, it is
well-known (\cite{5-cubic}) that the space $T^{(0,1)}_pM$ is
generated by the following holomorphic vector fields, tangent to $M$:
\[
\mathcal L_j:=\frac{\partial}{\partial \overline
z_j}+\sum_{l=1}^k\,\frac{\partial \overline\Xi_l}{\partial \overline
z_j}(z,\overline z,w)\,\frac{\partial}{\partial \overline w_l} \ \ \
\ \ {\scriptstyle (j=1,\ldots,n)}.
\]

\begin{Definition}
A CR-generic submanifold $M\subset\mathbb C^{n+k}$ of CR dimension $n$
and of codimension $k$ is called {\sl of finite type}
at a point $p \in M$ whenever the above
generators $\mathcal L_1,\ldots,\mathcal L_n,\overline{\mathcal
L}_1,\ldots,\overline{\mathcal L}_n$ together
with all of their Lie brackets of any length
span the complexified tangent space $\mathbb C
T_pM$ at the point $p$.
\end{Definition}

Of course, finite-typeness at a point is an open condition.
In \cite{Bloom}, Bloom and Graham introduced an effective method to
construct homogeneous CR manifolds that we now explain
briefly. Consider the complex space $\mathbb C^{n+k}$ equipped with
the variables $z_1,\ldots,z_n,w_1:=u_1+iv_1,\ldots,w_k:=u_k+iv_k$,
assign weight $1$ to the variables $z_i$ and assign some arbitrary
weights $\ell_j$ to each $w_j$ for $j=1,\ldots,k$. Then a CR manifold
$M\subset\mathbb C^{n+k}$ is called {\sl represented in Bloom-Graham
normal form} \cite{Bloom, Boggess} whenever it is defined as the
graph of some real-valued functions\footnote{There is also a more
general definition of Bloom-Graham normal form which one can find it
for example in \cite{Baracco}}:
\begin{equation}
\label{Bloom}
v_j
=
\Phi_j(z,\overline z,u_1,\ldots,u_{j-1})
+
{\rm o}(\ell_j) \ \ \ \ \ \ \
{\scriptstyle (j=1,\ldots,k),}
\end{equation}
where each $\Phi_j$ is a weighted homogeneous polynomial of the
weight $\ell_j$ enjoying the following two statements:

\begin{itemize}

\smallskip\item[$(i)$]
there are no {\it pure} terms $z^\alpha u^\beta$ or
$\overline z^\alpha u^\beta$ among the polynomials $\Phi_\bullet$ for
some integers $\alpha$ and $\beta$;

\smallskip\item[$(ii)$] for each $1\leqslant j<
i$ and for any nonnegative integers $\alpha_1,\ldots,\alpha_j$, the
polynomial $\Phi_i$ does not include any term of the form
$u_1^{\alpha_1}\cdots u_j^{\alpha_j}\Phi_j$.
\end{itemize}\smallskip

Every CR manifold represented in this form is of finite type ({\it
see} \cite{Boggess} page 181). Bloom and Graham also showed that
every CR manifold represented by the above expressions
\thetag{\ref{Bloom}} can be transformed to such a normal form by
means of some algebraic changes of coordinates ({\it see}
\cite{Bloom}, Theorem 6.2).

\begin{Definition}
A real analytic CR-generic
manifold $M\subset\mathbb C^N$ with coordinates $(Z_1, \dots,
Z_N)$ is called {\sl holomorphically nondegenerate} at $p\in M$ if
there is no local {\em nonzero} vector field of type $(1, 0)$:
\[
{\sf X}:=\sum_{j=1}^N\,f_j(Z_1,\ldots,Z_N) \frac{\partial}{\partial
Z_j}
\]
having coefficients $f_j$ holomorphic in
a neighborhood of $p$ such that ${\sf
X}|_M$ is tangent to $M$ near $p$.
\end{Definition}

Every {\em connected} real analytic generic CR manifold is either
holomorphically nondegenerate at every point or at no point
(\cite{ Merker-Porten-2006}).

\begin{Definition}
(\cite{AMS, Baouendi, Beloshapka1997, MS})
A {\sl (local)
infinitesimal CR-automorphism of $M$}, when understood extrinsically,
is a local holomorphic vector field:
\begin{equation}
\label{XX}
{\sf X}
=
\sum_{i=1}^n\,Z^i(z,w) \frac{\partial}{\partial
z_i} + \sum_{j=1}^k\,W^j(z,w) \frac{\partial}{\partial w_j}
\end{equation}
whose real part ${\rm Re}\, X = \frac{ 1}{ 2} ( {\sf X} +
\overline{\sf X} )$ is tangent to $M$.
\end{Definition}

The collection of all infinitesimal CR-automorphisms of $M$
constitutes a Lie
algebra which is called the {\sl Lie algebra of infinitesimal
CR-automorphisms} of $M$ and is denoted by $\frak{aut}_{CR}(M)$.

The
notion of holomorphically nondegeneracy was raised by Nancy Stanton in
\cite{Stanton} where she proved that
for a {\em hypersurface} $M \subset \mathbb{ C}^{n+1}$ (always generic),
$\mathfrak{ aut}_{ CR} ( M)$ is finite-dimensional if and only
of $M$ is holomorphically nondegenerate.
Amazingly enough, one
realizes that the concept of tangent vector
fields completely
independent of $\overline{ Z}_1, \dots, \overline{ Z}_N$
which points out a strong degeneracy
can in fact be traced
back at least to Sophus Lie's works ({\em cf.} pp.~13--14
of~\cite{ Merker-2010}).  In general codimension $k \geqslant 1$, the
Lie algebra $\mathfrak{ aut}_{ CR} ( M)$
of infinitesimal CR-automorphisms of
a CR-generic real analytic
$M \subset \mathbb{ C}^{ n+k}$ is finite-dimensional if and
only if $M$ is holomorphically nondegenerate {\em and} of finite type
(\cite{ Gaussier-Merker}).

Determining such Lie algebras $\mathfrak{ aut}_{ CR} ( M)$
is the same as knowing the {\em
CR-symmetries} of $M$, a question which lies at the heart of the
(open) problem of classifying all local analytic CR manifolds up to
biholomorphisms. In the groundbreaking works of Sophus Lie and
his followers (Friedrich Engel, Georg Scheffers, Gerhard Kowalewski,
Ugo Amaldi and others), the most fundamental question in
concern was to draw up lists of possible Lie algebras
$\mathfrak{aut}_{CR}(M)$ which would classify all possible $M$'s
according to their CR symmetries.


\subsection{The main idea}
\label{SHAM}

Now, let us explain briefly the main strategy used in \cite{HAMS} for
computing the Lie algebra $\frak{aut}_{CR}(M)$ associated to an
arbitrary real analytic generic CR manifold $M\subset\mathbb
C^{n+k}$, represented as the graph of the $k$ complex defining
equations as \thetag{\ref{Xi}}. For this, we shall proceed to do
successively the following steps:

\begin{itemize}
\item[$\bullet$] According to the definition, a holomorphic vector
field:
\[
{\sf X}=\sum_{j=1}^n\,Z^j(z,w)\,\partial
_{z_j}+
\sum_{l=1}^k\,W^l(z,w)\,\partial_{w_l}
\]
belongs to $\frak{aut}_{CR}(M)$ whenever it enjoys the tangency equations:
\begin{equation}
\label{Tangency} \footnotesize \aligned 0 \equiv (\overline{\sf
X}+{\sf X}) & \big[ \overline{w}_j+w_j -
\overline{\Xi}_j(\overline{z},z,w) \big] =
\\
& = \overline{\sf X}\,\big[\overline{w}_j+w_j -
\overline{\Xi}_j(\overline{z},z,w) \big] + {\sf X}\,
\big[\overline{w}_j+w_j - \overline{\Xi}_j(\overline{z},z,w) \big]
\\
& = \overline{W}^j(\overline{z},\overline{w}) -
\sum_{i=1}^n\,\overline{Z}^i(\overline{z},\overline{w})\,
\frac{\partial\overline{\Xi}_j}{\partial \overline{z}_i}
(\overline{z},z,w) +
\\
& \ \ \ \ \ + W^j(z,w) - \sum_{i=1}^n\, Z^i(z,w)\,
\frac{\partial\overline{\Xi}_j}{\partial z_i}(\overline{z},z,w) -
\sum_{l=1}^k\,W^l(z,w)\, \frac{\partial\overline{\Xi}_j}{\partial
w_l}(\overline{z},z,w)
\\
& \ \ \ \ \ \ \ \ \ \ \ \ \ \ \ \ \ \ \ \ \ \ \ \ \ \ \ \ \ \ \ \ \ \
\ \ \ \ \ {\scriptstyle{(j\,=\,1\,\cdots\,k).}}
\endaligned
\end{equation}
For each $j=1,\ldots,k$, let us refer to the above equality as the
{\sl $j$-th tangency equation} of $M$. Now, the Taylor series
formulaes:
\begin{equation}
\label{Taylor-1} \aligned Z^i(z,w) = \sum_{\alpha\in\mathbb
N^n}\,z^\alpha\,Z^{i}_\alpha(w) \ \ \ \ \ \text{\rm and} \ \ \ \ \
W^l(z,w) = \sum_{\alpha\in\mathbb N^n}\,z^\alpha\,W^{l}_\alpha(w),
\endaligned
\end{equation}
bring the tangency equation into the form:
\begin{equation}
\label{modified-tangency} \footnotesize \aligned 0 & \equiv
\sum_{\alpha\in\mathbb N^n}\,\overline{z}^\alpha\,
\overline{W}^{j}_\alpha \big(-w+\overline{\Xi}\big) -
\sum_{k=1}^n\,\sum_{\alpha\in\mathbb N^n}\, \overline{z}^\alpha\,
\overline{Z}^{k}_\alpha \big(-w+\overline{\Xi}\big)\,
\frac{\partial\overline{\Xi}_j}{\partial\overline{z}_k}
(\overline{z},z,w) +
\\
& \ \ \ \ \ + \sum_{\beta\in\mathbb N^n}\, z^\beta\,W^{j}_\beta(w) -
\sum_{k=1}^n\,\sum_{\beta\in\mathbb N^n}\, z^\beta\,Z^{k}_\beta(w)\,
\frac{\partial\overline{\Xi}_j}{\partial z_k} (\overline{z},z,w) -
\sum_{l=1}^k\,\sum_{\beta\in\mathbb N^n}\, z^\beta\,W^{l}_\beta(w)\,
\frac{\partial\overline{\Xi}_j}{\partial w_l} (\overline{z},z,w)
\\
& \ \ \ \ \ \ \ \ \ \ \ \ \ \ \ \ \ \ \ \ \ \ \ \ \ \ \ \ \ \ \ \ \ \
\ \ \ \ \ \ \ \ \ \ \ \ \ \ \ \ \ \ \ \ \ \ \ \ \ \
{\scriptstyle{(j\,=\,1\,\cdots\,k)}}.
\endaligned
\end{equation}
\item[$\bullet$] After this, one
replaces the appearing functions $\overline
W^{\bullet}_\bullet(-w+\overline{\Xi})$ and $\overline
Z^{\bullet}_\bullet(-w+\overline{\Xi})$ according to the following
slightly artificial expansions:
\[
\small \aligned \overline{A}\big(-w+\overline{\Xi}\big) & =
\overline{A}\big(w+(-2w+\overline{\Xi})\big)
\\
& = \sum_{\gamma\in\mathbb N^d}\,
\frac{\partial^{\vert\gamma\vert}\overline{A}}{
\partial w_\gamma}(w)\,
\frac{1}{\gamma!}\,
\big(-2w+\overline{\Xi}(\overline{z},z,w)\big)^\gamma,
\endaligned
\]
and one next substitutes each $-2w_j + \overline{ \Xi}_j$ with:
\[
-2w_j + \overline{\Xi}_j(\overline{z},z,w) = \sum_{\alpha\in\mathbb
N^n}\,\sum_{\beta\in\mathbb N^n}\, \overline{z}^\alpha\,z^\beta\,
\overline{\Xi}_{j,\alpha,\beta}(w) \ \ \ \ \ \ \ \ \ \ \ \ \
{\scriptstyle{(j\,=\,1\,\cdots\,k)}}.
\]
\item[$\bullet$]
Modifying the equations \thetag{\ref{modified-tangency}} by the two
already presented formulaes, we reach $k$ homogenous equations so
that their right hand sides are some certain combinations of the
(yet unknown) functions $Z^i_\alpha(w)$ and $W^l_\beta(w)$, their
derivations and also the variables $w_1,\ldots,w_k$ with the
coefficients of the form $z^\mu\overline{z}^\nu$ for
$\mu=(\mu_1,\ldots,\mu_n)$ and $\nu=(\nu_1,\ldots,\nu_n)$ (here we
set $z_\mu:=z_1^{\mu_1}\ldots z_n^{\mu_n}$). To satisfy these
equations, one should extract the coefficients of the various
monomials $z^\mu\overline{z}^\nu$ and equate them to zero. Then, one
attains a (usually lengthy and complicated) homogeneous linear system
of complex partial differential equations with the unknowns
$Z^{i}_\alpha(z,w)$ and $W^{l}_\alpha(z,w)$\,\,---\,\,namely a linear
system of differential polynomials of the differential algebra
$R:=\mathbb C(w,\overline w)[Z^i_\alpha,W^l_\beta]$. The solution of
this system yields the desired coefficients $Z^i(z,w)$ and $W^l(z,w)$
in \thetag{\ref{Taylor-1}}. \item[$\bullet$] To solve this already
constructed {\sc pde} system, we employ the effective tools of {\it
differential algebra}, equipped with some additional operators such
as {\it bar-reduction}. For more details we refer the reader to
\cite{HAMS}.
\end{itemize}

Let us conclude this section by a short discussion of the abilities
and advantages of the algorithm introduced in \cite{HAMS} and also
the difficulties and weaknesses inherent in it. This algorithm
manages to compute {\it systematically} the
desired algebras associated to {\it arbitrary} CR manifolds
$M\subset\mathbb C^{n+k}$. Moreover, it employs the powerful
techniques of differential algebra and the ability of computer
algebra to provide a more effective method. In fact, it handles more
appropriately the most complicated part of the computations, namely
solving the associated {\sc pde} systems ({\it cf.}
\cite{Beloshapka1997,Mamai2009,5-cubic,MS,Shananina2000}).
Specifically, for the significant class of
{\it rigid} CR manifolds\,\,---\,\,those whose defining
equations $\Xi$ are independent of the variables $w$\,\,---\,\,
the computations are considerabl eased up. {\it
However}, the main obstacle we encounter is that, because of
the complexity of the {\sc pde} systems,
as much as the number of variables $z_i$ and $w_l$
increases, then the cost of the associated computation grows
extensively and the implementation of the algorithm rapidly reveals
limits concerning the capacity of computer systems.

\section{Rigidity of homogeneous CR manifolds and
\\
similarity in the expressions of the corresponding sought functions}
\label{section rigidity-similarity}

Consider the complex space $\mathbb C^{n+k}$ equipped with $n$
complex variables $z_1,\ldots,z_n$ of weight one, identically, and
$k$ complex variables $w_1:=u_1+iv_1,\ldots,w_k:=u_k+iv_k$ of certain
weights $1<[w_1]\leqslant [w_2]\leqslant \ldots\leqslant  [w_k]$ and define the
homogeneous manifold $M$  of CR dimension $n$ and codimension $k$ as:
\begin{equation}
\aligned \label{model-model} M:=\left\{(z,w): \ \ \left[
\begin{array}{l}
\Xi_1(v_1,z,\overline z,u):=v_1-\Phi_1(z,\overline z,u)\equiv 0,
\\
\Xi_2(v_2,z,\overline z,u):=v_2-\Phi_2(z,\overline z,u)\equiv 0,
\\
\ \ \ \ \ \ \ \vdots
\\
\Xi_j(v_j,z,\overline z,u):=v_j-\Phi_j(z,\overline z,u)\equiv 0,
\\
\ \ \ \ \ \ \ \vdots
\\
\Xi_k(v_k,z,\overline z,u):=v_k-\Phi_k(z,\overline z,u)\equiv 0,
\end{array}
\right.\right\},
\endaligned
\end{equation}
with the weighted homogeneous polynomials $\Phi_j$ of the certain
weight $[w_j]$ for $j=1,\ldots,k$. {\it Throughout this paper} we
assume that these manifolds are holomorphically nondegenerate and of
finite type, which guarantees the finite dimensionality of their
associated algebras of infinitesimal CR-automorphisms.

\subsection{Gradation and polynomiality}

At first, let us show two significant intrinsic features of the
already mentioned homogeneous CR manifolds, namely gradation (in the
sense of Tanaka) and polynomiality of their associated algebras of
infinitesimal CR-automorphisms.

For a CR manifold $M$ as above, consider the holomorphic vector
field:
\[
{\sf X}:=\sum_{j=1}^n Z^j(z,w)\partial_{z_j}+\sum_{l=1}^k
W^l(z,w)\partial_{w_l}
\]
as an element of $\frak g:=\frak{aut}_{CR}(M)$. Since the above
coefficients $Z^j$ and $W^l$ are all holomorphic, then one can expand
them as their Taylor series and thus decompose $\sf X$ into its
weighted homogeneous components as follows:
\begin{equation}
\label{X-decomppse} {\sf X}:={\sf X}_{-\rho}+\cdots+{\sf X}_{-1}+{\sf
X}_0+{\sf X}_1+\cdots+{\sf X}_t+\cdots, \ \ \ \ \ {\scriptstyle
(\rho\,,\, t\,\in\,\mathbb N)}.
\end{equation}
We need the following facts:

\begin{Lemma}
The minimum homogeneous component ${\sf X}_{-\rho}$ in the above
decomposition of $\sf X$ is of the weight $-\rho=-[w_k]$, where $w_k$
has the maximum homogeneity among the complex variables appearing in
\thetag{\ref{model-model}}.
\end{Lemma}

\proof It is just enough to observe that the tangent space of
holomorphic vector fields can be generated by the standard fields
$\partial_{z_j},\partial_{w_l}$ of the certain weights $-[z_j]$ and
$-[w_l]$ for $j=1,\ldots,n$ and $l=1,\ldots,k$. Among these standard
fields, the minimum homogeneity belongs to $\frac{\partial}{\partial
w_k}$. This completes the proof.
\endproof

\begin{Lemma}
\label{lemma-polynomiality} Each of the above weighted homogeneous
components ${\sf X}_t, \, t\geqslant -\rho$ is again an infinitesimal
CR-automorphism, namely belongs to $\frak g$.
\end{Lemma}

\proof Since $\sf X$ is an infinitesimal CR-automorphism then we
have:
\[ 0\equiv \big({\sf X}+\overline{\sf X}\big)|_{\Xi_j},
\ \ \ \ {\scriptstyle (j=1,\ldots,k)}.
\]
Now, for each component ${\sf X}_t$ of homogeneity $t$ and since each
defining function $\Xi_j$ is homogeneous of weight $\ell_j:=[w_j]$
then one verifies that the polynomial $({\sf X}_t+\overline{\sf
X}_t)|_{\Xi_j}$ is either zero or a homogeneous polynomial of weight
$t+\ell_j$. Hence we have:
\begin{equation}
\label{sum} \aligned 0\equiv \big({\sf X}+\overline{\sf
X}\big)|_{\Xi_j}=\underbrace{\big({\sf X}_{-\rho}+\overline{\sf
X}_{-\rho}\big)|_{\Xi_j}}_{P_{\rho,j}}+& \ldots+\underbrace{\big({\sf
X}_{t}+\overline{\sf X}_{t}\big)|_{\Xi_j}}_{P_{t,j}}+\cdots \ \ \ \ \
\ \ \ {\scriptstyle (j=1,\ldots,k)},
\endaligned
\end{equation}
in which each polynomial function $P_{t,j}$ is either zero or a
homogeneous polynomial of the weight $t+\ell_j$. Hence, we have some
certain weighted homogeneous polynomials $P_\bullet$ with distinct
homogeneities and consequently with linear independency. Hence, one
obtains from \thetag{\ref{sum}} that:
\[
0\equiv P_{t,j}(z,w)=\big({\sf X}_{t}+\overline{\sf
X}_{t}\big)|_{\Xi_j}, \ \ \ {\scriptstyle t\geqslant -\rho, \ \ \
(j=1,\ldots,k)}.
\]
This equivalently means that each component ${\sf X}_t$ of $\sf X$ is
an infinitesimal CR-automorphism.
\endproof

Now, we can prove the polynomiality of the sought algebras:

\begin{Corollary}
\label{Polynomiality} If:
\[
{\sf X}_t:=\sum_{j=1}^n Z_t^j(z,w)\partial_{z_j}+\sum_{l=1}^k
W_t^l(z,w)\partial_{w_l}
\]
 is a weighted homogeneous CR-automorphism of
weight $t\geqslant  -\rho$ then, its coefficients $Z_t^j(z,w)$ and
$W_t^l(z,w)$ are weighted homogeneous polynomials of the weights
$0,\ldots,t-1$.
\end{Corollary}

\proof It is a straightforward consequence of the decomposition
\thetag{\ref{X-decomppse}}.
\endproof

Now, let us consider the gradation of $\frak g$. First we need the
following definition:

\begin{Definition}
\label{rigidity}
 The Lie algebra $\frak g:=\frak{aut}_{CR}(M)$ of an arbitrary
CR manifold is called {\sl graded} in the sense of Tanaka, whenever
it can be expressed in the form:
\[
\frak g:=\frak g_{-\rho}\oplus\frak
g_{-\rho+1}\oplus\cdots\oplus\frak g_{-1}\oplus\frak g_0\oplus\frak
g_1\oplus\cdots\oplus\frak g_{\varrho}, \ \ \ \ \ \ \ \ \
{\scriptstyle \rho,\,\varrho\,\in\,\mathbb N}
\]
with $[\frak g_i,\frak g_j]\subset\frak g_{i+j}$. Furthermore, we say
that $M$ has {\sl rigidity} if the positive part:
\[
\frak g_+:=\frak g_1\oplus\ldots\oplus\frak g_{\varrho}
\]
of $\frak{aut}_{CR}(M)$ is zero.
\end{Definition}

\begin{Proposition}
\label{Gradation} For a finite type holomorphically nondegenerate CR manifold
$M\subset\mathbb C^{n+k}$ represented by the above defining function
\thetag{\ref{model-model}}, the associated Lie algebra
$\frak g$ of infinitesimal CR-automorphisms is a graded Lie algebra,
in the sense of Tanaka, of the form:
\[
\frak g:=\frak g_{-\rho}\oplus\ldots\oplus\frak g_{-1}\oplus\frak
g_0\oplus\frak g_1\oplus\ldots\oplus\frak g_\varrho, \ \ \ \
{\scriptstyle \rho\,,\,\varrho\,\in\,\mathbb N}.
\]
\end{Proposition}

\proof According to the above two lemmas, if $\frak g_t$ is the
collection of infinitesimal CR-automorphisms of weight $t$, then
$\frak g$ admits a gradation like:
\[
\frak g:=\frak g_{-\alpha}\oplus\ldots\oplus\frak g_{-1}\oplus\frak
g_0\oplus\frak g_1\oplus\cdots\oplus\frak g_t\oplus\cdots.
\]
Furthermore, holomorphically nondegeneracy and finite typness of $M$
guarantees that $\frak g$ is finite dimensional and hence there
exists an integer $\varrho$ in which $g_\beta\equiv 0$ for each
$\beta>\varrho$. Now, it remains to show that this gradation is in
the sense of Tanaka. Namely, we shall prove that for each two
homogeneous infinitesimal CR-automorphisms ${\sf X}_1$ and ${\sf
X}_2$ of certain homogeneities $t_1$ and $t_2$, the Lie bracket
$[{\sf X}_1,{\sf X}_2]$ belongs to $\frak g_{t_1+t_2}$. For this, it
is enough to show this statement for vector fields of the forms ${\sf
X}_i:=F_i(z,w)\frac{\partial}{\partial x_i}, i=1,2$ where $x_i$ is a
complex variable $z_\bullet$ or $w_\bullet$. According to the
homogeneities of these vector fields, $F_1$ and $F_2$ are two
polynomials of the weights $t_1+[x_1]$ and $t_2+[x_2]$, respectively.
Now, we have:
\[
[{\sf X}_1,{\sf X}_2]=F_1\,\frac{\partial F_2}{\partial
x_1}\frac{\partial}{\partial x_2}-F_2\,\frac{\partial F_1}{\partial
x_2}\frac{\partial}{\partial x_1}.
\]
In this expression, the derivations $\frac{\partial F_2}{\partial
x_1}$ and $\frac{\partial F_1}{\partial x_2}$ are either zero or
homogeneous polynomials of the weights $t_2+[x_2]-[x_1]$ and
$t_2+[x_1]-[x_2]$, respectively. Now, simple simplifications yields
that the above Lie bracket is a homogeneous vector field of the
weight $t_1+t_2$, as desired.
\endproof

\subsection{A comparison between the so far achieved results}

For a CR manifold $M\subset\mathbb C^{n+k}$, of {\it fixed}
CR-dimension $n$ and codimension $k$ represented in coordinates
$z_1,\ldots,z_n,w_1,\ldots,w_k$ as \thetag{\ref{model-model}},
consider the infinitesimal CR-automorphism:
\begin{equation*}
{\sf
X}:=\sum_{j=1}^n\,Z^j(z,w)\partial_{z_j}+\sum_{l=1}^k\,W^l(z,w)\partial_{w_l}.
\end{equation*}
In this subsection, we have a general look at the results through
literatures. In \cite{Shananina2000}, Shananina computed the Lie
algebras of infinitesimal CR-automorphisms associated to the
Beloshapka's homogeneous models\,\,---\,\,included in the class of
under consideration CR manifolds\,\,---\,\,with $n=1$ and
$k=3,\ldots,7$ and presented the expressions of the associated
holomorphic coefficients $Z(z,w)$ and $W^l(z,w)$, $l=1,\ldots,k$ in
Theorem 1 page 386 of his paper (he used the notations $f(z,w)$ and
$g^l(z,w)$, instead). Upon a close scrutiny, one observes a nice {\it
similarity} in the expressions of those computed functions. More
precisely, all the expressions of $Z(z,w)$ and $W^l(z,w)$, presented
in the various types $k=3$ to $k=7$ are nearly the same. For
instance, in the two cases $k=3$ and $k=4$ represented by the
following defining equations:
\begin{equation*}
\aligned {\tiny\big(\begin{array}{l} n=1
\\
 k=3
\end{array}
\big) }
 \ :\left\{
\begin{array}{l}
v_1=z\overline z,
\\
 v_2=z^2\overline z+z\overline z^2,
\\
v_3=-i(z^2\overline z-z\overline z^2),
\end{array}
\right. \ \ \ \ \ \ \ \ \ {\tiny\big(\begin{array}{l} n=1
\\
 k=4
\end{array}
\big): }\ \left\{
\begin{array}{l}
v_1=z\overline z,
\\
 v_2=z^2\overline z+z\overline z^2,
\\
v_3=-i(z^2\overline z-z\overline z^2),
\\
v_4=z^3\overline z+z\overline z^3+\lambda\,z^2\overline z^2, \ \ \
{\scriptstyle \lambda\in\mathbb C}
\end{array}
\right.
\endaligned
\end{equation*}
 one obtains the
expressions:
\begin{equation}
\footnotesize
\label{k=3,4}
 \aligned \boxed{n=1,\,k=3}&:
\\
 &Z(z,w)={\sf a}+i\,{\sf b}+({\sf d}+i\,{\sf
e})\,z, \ \ \ \ \ W^1(z,w)={\sf c}_1+2\,({\sf b}+i\,{\sf
a})\,z+2\,{\sf d}w_1,
\\
& W^2(z,w)={\sf c}_2+2\,({\sf b}+i\,{\sf a})\,z^2+4\,{\sf
a}w_1+3\,{\sf d}w_2-{\sf e}w_3,
\\
& W^3(z,w)={\sf c}_3+2\,({\sf a}-i\,{\sf b})\,z^2+4\,{\sf b}w_1+{\sf
e}w_2+3\,{\sf d}w_3,
\\
\boxed{n=1,\,k=4}&:
\\
 &Z(z,w)={\sf a}+i\,{\sf b}+{\sf d}z, \ \ \ \ \ W^1(z,w)={\sf c}_1+2\,({\sf b}+i\,{\sf
a})\,z+2\,{\sf d}w_1,
\\
& W^2(z,w)={\sf c}_2+2\,({\sf b}+i\,{\sf a})\,z^2+4\,{\sf
a}w_1+3\,{\sf d}w_2,
\\
& W^3(z,w)={\sf c}_3+2\,({\sf a}-i\,{\sf b})\,z^2+4\,{\sf
b}w_1+3\,{\sf d}w_3,
\\
& W^4(z,w)={\sf c}_4+2\,({\sf b}+i\,{\sf a})\,z^3+(2\lambda+3)\,{\sf
a}w_2+2\,(\lambda-3)\,{\sf b}w_3+4\,{\sf d}w_4,
\\
& \ \ \ \ \ \ \ \ \ \ \ \ \ \ \ \ \ \ \ \ \ \ \ \ \ \ \ \ \ \ \ \ \ \
\ \ \ \ \ \ {\scriptstyle ({\sf a,b,c}_i,{\sf d,e}\in\mathbb R)}.
\endaligned
\end{equation}
A plain comparison between the above two sets of the functions
$Z,W^1,W^2,W^3$ can demonstrate the mentioned similarity. For
example, even if a function is not equal to its correspondence,
the degrees of each term in both of them and the variables which
appear inside are same. This type of similarity and harmony is
visible among all the expressions obtained in \cite{Shananina2000}
for $n=1$ and $k=3,\ldots,7$.
Moreover, a glance on the
Beloshapka's results for his model with $n=1,\,k=2$ in
\cite{Beloshapka1997}, page 402 ({\it see} the expressions of $f, g$
and $h$ which are $Z,W^1$ and $W^2$ in our notations) indicates such
similarity, too. In contrast, the computations of the infinitesimal
CR-automorphisms associated to the CR manifold with $n,k=1$,
presented in \cite{MS}, says that the expressions of the two
functions $Z(z,w)$ and $W^1(z,w)$ are fairly different. Namely we
have:
\begin{equation*}
\footnotesize\label{(1.1)}\aligned Z(z,w_1)&={\sf a}+i{\sf b}+({\sf
c}+i{\sf d})\,w_1+ (\frac{1}{2}\,{\sf f}+i\,{\sf e}+{\sf
h}\,w_1)\,z+(2\,{\sf d}+2\,i\,{\sf c})\,z^2,
\\
W^1(z,w_1)&={\sf f}+{\sf g}\,w_1+{\sf h}\,w_1^2+\big(2\,i\,({\sf
a}-i{\sf b})+2\,i\,({\sf c}-i{\sf d})\,w_1\big)\,z,
\endaligned
\end{equation*}
for some eight real integers ${\sf a},{\sf b},{\sf c},\ldots,{\sf
h}$. In more details, in the above expression of $Z(z,w_1)$ we see
the variable $w_1$ while the functions $Z(z,w)$ in the next types are
independent of the variables $w_i$. Moreover, the degree of the
variable $z$ in $Z(z,w)$ in all already presented cases $n=1$ and
$k=2,\ldots,7$ is 1 while in this case it is of degree 2 (compare the
above expression of $Z(z,w_1)$ with \thetag{\ref{k=3,4}}).

On the other hand\,\,---\,\,and {\it parallel} to the similarity of
the obtained holomorphic functions\,\,---\,\,we know that all the
mentioned CR manifolds with $n=1$ and $k=2,\ldots,7$ have {\it
rigidity} while for $n,k=1$ it does not retain this property
\cite{Mamai2009}.

We {\it claim} that the similarity in the expressions of the
functions $Z^i(z,w)$ and $W^l(z,w)$ associated to the CR manifolds of
CR-dimension and codimensions $n=1$ and $k=2,\ldots,7$ has a close
relationship to their rigidity property. Before considering this
assertion, let us illustrate by the following example that how one
can extract the infinitesimal CR-automorphisms $\sf X$ from the
obtained expressions of $Z(z,w)$ and $W^l(z,w)$.

\begin{Example}
\label{(1,3)} According to \cite{5-cubic} ({\it see} also the case
$K=3$, \cite{Shananina2000}, Theorem 1 for another presentation), an
infinitesimal CR-automorphism for the CR manifold $M\subset\mathbb
C^{1+3}$, represented in coordinates $\mathbb C\{z,w_1,w_2,w_3\}$ as
the graph of the defining functions:
\[
\left\{ \aligned w_1-\overline{w}_1 & = 2i\,z\overline{z}
\\
w_2-\overline{w}_2 & = 2i\,z\overline{z}(z+\overline{z})
\\
w_3-\overline{w}_3 & = 2\,z\overline{z}(z-\overline{z})
\endaligned\right.
\]
 is a holomorphic vector field:
\[
{\sf X}:=Z(z,w)\,\partial_ z+\sum_{l=1}^3\,W^l(z,w)\,\partial_{w_l},
\]
with the desired coefficients ({\it cf.} \thetag{\ref{k=3,4}})
\begin{equation*}
\footnotesize\aligned Z(z,w)&={\sf a}+i\,{\sf b}+({\sf d}+i\,{\sf
e})\,z,
\\
W^1(z,w)&={\sf c}_1+2\,({\sf b}+i\,{\sf a})\,z+2\,{\sf d}\,w_1,
\\
W^2(z,w)&={\sf c}_2+2\,({\sf b}+i\,{\sf a})\,z^2+4\,{\sf
a}\,w_1+3\,{\sf d}\,w_2-{\sf e}\,w_3,
\\
W^3(z,w)&={\sf c}_3+2\,({\sf a}-i\,{\sf b})\,z^2+4\,{\sf b}\,w_1+{\sf
e}\,w_2+3\,{\sf d}\,w_3,
\endaligned
\end{equation*}
for some seven real integers ${\sf a,b,c}_1,{\sf c}_2,{\sf c}_3,{\sf
d,e}$. Throughout this paper, we call such integers by the {\sl free
parameters}. Then, the Lie algebra $\frak{aut}_{CR}(M)$ is seven
dimensional, with the basis elements extracted as the coefficients of
the seven free parameters in the above general form of $\sf X$.
Namely, it is generated by the following seven holomorphic vector
fields:
\begin{equation}
\label{K=3}\footnotesize\aligned \boxed{\sf a}: \ \ \ \ \ \ \ {\sf
X}_1&:=\partial_z+2\,i\,z\partial_{w_1}+(2\,i\,z^2+4\,w_1)\,\partial_{w_2}+2\,z^2\,\partial_{w_3},
\\
\boxed{\sf b}: \ \ \ \ \ \ \ {\sf
X}_2&:=i\,\partial_z+2\,z\,\partial_{w_1}+2\,z^2\,\partial_{w_2}+(-2\,i\,z^2+4\,w_1)\,\partial_{w_3},
\\
\boxed{{\sf c}_1}: \ \ \ \ \ \ \ {\sf X}_3&:=\partial_{w_1},
\\
\boxed{{\sf c}_2}: \ \ \ \ \ \ \ {\sf X}_4&:=\partial_{w_2},
\\
\boxed{{\sf c}_3}: \ \ \ \ \ \ \ {\sf X}_5&:=\partial_{w_3},
\\
\boxed{\sf d}: \ \ \ \ \ \ \ {\sf
X}_6&:=z\,\partial_z+2\,w_1\,\partial_{w_1}+3\,w_2\,\partial_{w_2}+3\,w_3\,\partial_{w_3},
\\
\boxed{\sf e}: \ \ \ \ \ \ \ {\sf
X}_7&:=i\,z\,\partial_z-w_3\,\partial_{w_2}+w_2\,\partial_{w_3}.
\endaligned
\end{equation}
 The weights
associated in \cite{Shananina2000} to the appearing complex variables
are:
\[
[z]=1, \ \ \ [w_1]=2, \ \ \  [w_2]=[w_3]=3.
\]
Hence, a close look at the obtained basis holomorphic vector fields
${\sf X}_1,\ldots,{\sf X}_7$ gives the following weighted
homogeneities for them:

\medskip
\begin{center}
\begin{tabular} [t] { l | l l l l l l l}
 $\sf X$ & ${\sf X}_1$ & ${\sf X}_2$ & ${\sf X}_3$ & ${\sf X}_4$ & ${\sf
 X}_5$ & ${\sf X}_6$ & ${\sf X}_7$
\\
\hline {\textrm Hom} & -1 & -1 & -2 & -3 & -3 & 0 & 0
\end{tabular}
\end{center}
Therefore, the Lie algebra $\frak{aut}_{CR}(M)$ can be represented
as:
\[
\frak{aut}_{CR}(M):=\frak g_{-3}\oplus\frak g_{-2}\oplus\frak
g_{-1}\oplus\frak g_0,
\]
with $\frak g_{-3}:=\langle {\sf X}_4,{\sf X}_5\rangle$, with $\frak
g_{-2}:=\langle {\sf X}_3\rangle$, with $\frak g_{-1}:=\langle {\sf
X}_1,{\sf X}_2\rangle$ and with $\frak g_0:=\langle {\sf X}_6,{\sf
X}_7\rangle$. Now, one easily verifies that $M$ is rigid.
\end{Example}

\subsection{Homogeneous components}

Now, let us inspect the structure of homogeneous components of the
desired algebras of infinitesimal CR-automorphisms. At first, the
following proposition demonstrates the influence of the rigidity of
the under consideration homogeneous CR manifold $M$ on the structure.
It also reveals the close connection between the similarity property
of the holomorphic vector fields $Z^i$ and $W^l$ of  the CR manifolds
discussed in the last subsection and their rigidity.

\begin{Proposition}
\label{Rigidity-test} Let $M$ be a holomorphically nondegenerate
CR manifold of CR-dimension $n$ and codimension $k$ represented as
\thetag{\ref{model-model}}. Then, $M$ has rigidity if and only if for
any weighted homogeneous infinitesimal CR-automorphism:
\[
{\sf
X}:=\sum_{j=1}^n\,Z^j(z,w)\,\partial_{z_j}+\sum_{l=1}^k\,W^l(z,w)\,\partial_{w_l}
\]
of $M$, each weighted homogeneous polynomial $Z^j(z,w)$ (respectively
$W^l(z,w)$) is of weight at most $[z_i]$ (respectively of weight at
most $[w_l]$). In particular, $Z^j(z,w)$ (respectively $W^l(z,w)$) is
independent of the variables of the weights $\gvertneqq [z_j]$
(respectively of the weights $\gvertneqq [w_l]$).
\end{Proposition}

\proof If $\sf X$ is of weight homogeneity $d$, then since the
standard fields $\partial_{z_j}$ and $\partial_{w_l}$ have the
constant weights $-[z_i]$ and $-[w_l]$, respectively, then the
weighted homogeneous polynomials $Z^j$ and $W^l$ have the constant
degrees $d+[z_i]$ and $d+[w_l]$, respectively. Now, assume that $M$
has rigidity. The dimension of the Lie algebra $\frak{aut}_{CR}(M)$
is equal to the number of free parameters involved in the expressions
of the functions $Z^j(z,w), i=1,\ldots,n$ and $W^l(z,w),
l=1,\ldots,k$ ({\it see} Example \ref{(1,3)}). Each generator is
extracted from one of such free parameters as its coefficient in the
general form ${\sf
X}:=\sum_{j=1}^n\,Z^j(z,w)\,\partial_{z_j}+\sum_{l=1}^k\,W^l(z,w)\,\partial_{w_l}$.
In such expression, coefficients of the standard fields
$\partial_{z_j}$ come from the found functions $Z^j(z,w)$. Now,
rigidity of $M$ means to have no any (homogeneous) basis element
belonging to $\frak{aut}_{CR}(M)$ of the positive weighted
homogeneity. Hence, when we have standard field $\partial_{z_j}$ of
the weight $-[z_j]$, then no term of weight bigger than $[z_j]$
appears in its coefficient. Consequently, $Z^j(z,w)$\,\,---\,\,which
provides the coefficients of $\partial_{z_j}$ in the basis
elements\,\,---\,\,is independent of the variables of the weights
bigger than $[z_j]$. Similar fact holds when one considers the
coefficients $W^l(z,w)$ of the standard fields $\partial_{w_l}$.

\noindent For the converse, if none of the (weighted homogeneous)
holomorphic coefficients $Z^j(z,w), j=1,\ldots,n$ admits the terms of
weight larger than $[z_j]$, then the weight $[Z^j(z,w)]-[z_j]$ of
each term $Z^j(z,w)\partial_{z_j}$ of $\sf X$ is non-positive.
Similar fact holds for the terms $W^l(z,w)\partial_{w_l},
l=1,\ldots,k$. Consequently, $\frak{aut}_{CR}(M)$ does not contain
any (weighted homogeneous) basis element $\sf X$ of the positive
weight. In other words, $M$ has rigidity.
\endproof

For a graded Lie algebra:
\[
\frak g=\frak g_{-\rho}\oplus\cdots\oplus\frak g_{-1}\oplus\frak
g_0\oplus\frak g_1\oplus\cdots\oplus\frak g_\varrho,
\]
let us denote by $\frak g^{(t)}$ the graded subspace:
\[
\frak g^{(t)}:= \frak g_{-\rho}\oplus\cdots\oplus\frak g_{t}, \ \ \ \
{\scriptstyle (t=-\rho,\ldots,\varrho)}.
\]
According to definition of the graded algebras, one easily convinces
oneself that for $t=-\rho,\ldots,0$ each subspace $\frak g^{(t)}$ is
in fact a Lie subalgebra of $\frak g$.
 The idea behind the proof of Proposition
\ref{Rigidity-test} can lead one to obtain the following more general
conclusion.

\begin{Proposition}
\label{main-theorem} \label{g(t)} Let $M$ be a homogeneous
CR manifold of CR-dimension $n$ and codimension $k$ represented as
\thetag{\ref{model-model}}. Let $\frak g=\frak{aut}_{CR}(M)$ be of
the graded form:
\begin{equation}
\label{g}
 \frak g=\frak g_{-\rho}\oplus\cdots\oplus\frak
g_{-1}\oplus\frak g_0\oplus\frak g_1\oplus\cdots\oplus\frak
g_\varrho
\end{equation}
and let the (weighted homogeneous) infinitesimal CR-automorphism:
\begin{equation*} {\sf
X}=\sum_{j=1}^n\,Z^j(z,w)\,\partial_{z_j}+\sum_{l=1}^k\,W^l(z,w)\,\partial_{w_l}
\end{equation*}
belongs to $\frak g$. Then,
\begin{itemize}
\item[$(i)$] $\sf X$ belongs to $\frak g_t$ for
$t=-\rho,\ldots,\varrho$, if and only if each coefficient $Z^j(z,w)$
(respectively $W^l(z,w)$) is homogeneous of the precise weight
$[z_j]+t$ (respectively $[w_l]+t$). In particular, each $Z^j(z,w)$
(respectively $W^l(z,w)$) is independent of the variables of weights
$\gneqq [z_j]+t$ (respectively of weights $\gneqq [w_l]+t$).
\item[$(ii)$] $\sf X$ belongs to $\frak g^{(t)}$ if and only if each
coefficient $Z^j(z,w)$ (respectively $W^l(z,w)$) is homogeneous of
weight at most $[z_j]+t$ (respectively of weights at most $[w_l]+t$).
Specifically, each function $Z^j(z,w)$ (respectively, $W^l(z,w)$) is
independent of the variables of weights $\gneqq [z_j]+t$)
(respectively of weights $\gneqq [w_l]+t$)\item[$(iii)$] In
particular, the negative part $\frak g_-=\frak g^{(-1)}$ of the
CR manifold $M$ is generated by the elements $\sf X$ of $\frak g$
with the coefficients $Z^j(z,w)$ (respectively $W^l(z,w)$)
independent of the variables of weights $\gneqq [z_j]-1$
(respectively of weights $\gneqq [w_l]-1$).
\end{itemize}
\end{Proposition}

\proof The proof is similar to that of Proposition
\ref{Rigidity-test}. Here, we prove the first item (i). If a
holomorphic vector field $\sf X$ belongs to $\frak g_t$, then it is
homogeneous of the weight $t$. In the expression of $\sf X$, each
standard field $\partial_{z_j}$ is of the fixed weight $-[z_j]$ and
hence having the field $Z^j(z,w)\,\partial_{z_j}$ of the precise
homogeneity $t$, then the coefficient $Z^j(z,w)$ must be of the weigh
$[z_j]+t$. Similar conclusion holds for the functions $W^l(z,w)$. The
converse of the assertion can be concluded in a very similar way. In
particular, since all the variables in $Z^j(z,w)$ have the positive
weight then, each function $Z^j(z,w)$\,\,---\,\,which must be of the
weight $[z_j]+t$\,\,---\,\,is independent of the variables of the
weights bigger than $[z_j]+t$.
\endproof

\begin{Remark}
\label{Remark-algorithm} This proposition suggests an appropriate way
to compute each subspace $\frak g^{(t)}$ (specifically, to compute
$\frak g^{(0)}=\frak g_-\oplus\frak g_0$ in the rigid case). Namely,
it states that to find each $\frak g^{(t)}$ it is not necessary to
compute the surrounding algebra $\frak g=\frak{aut}_{CR}(M)$ of
infinitesimal CR-automorphisms, but it is sufficient to compute just
the (homogeneous) coefficients $Z^j(z,w)$ and $W^l(z,w)$ of
homogeneities less than $[z_j]+t$ and $[w_l]+t$, respectively. This
extremely reduces the size of computation. Similar process can be
employed when we aim to compute only the $t$-th component $\frak g_t$
of $\frak g$. Moreover, one can divide the general computation of the
Lie algebra $\frak g$ of infinitesimal CR-automorphisms into some
distinct sub-computations of its components $\frak g_t$ for
$t=-\rho,\ldots,\varrho$. In the next section, we will use the
results to provide an algorithm for computing $\frak{aut}_{CR}(M)$.
\end{Remark}

\section{Computing the homogeneous components}

 \label{section computation the components}

As mentioned in Remark \ref{Remark-algorithm}, to compute the
holomorphic coefficients $Z^j(z,w)$ and $W^l(z,w)$ of the vector
fields ${\sf X}\in\frak g^{(t)}$, one can assume these functions
independent of the variables with the associated weights larger than
$[z_j]+t$ and $[w_l]+t$, respectively. In this section, we aim to
develop this result for constructing a very convenient method of
computing each subspace $\frak g^{(t)}$ and each component $\frak
g_t$ associated to our homogeneous CR manifold.

\subsection{Computing each component $\frak g_t$}
\label{strategy-1}

For each element $\sf X$ of the $t$-th component $\frak g_t$,
Proposition \ref{main-theorem} enables one to attain an upper bound
for the weight degree of each of its desired (polynomial)
coefficients $Z^j(z,w)$ and $W^l(z,w)$. Hence, we can predict the
expression of these polynomials as the elements of the polynomial
ring $\mathbb C[z,w]$. Then finding these expressions, it is
necessary and sufficient to seek their constant coefficients. One can
pick the following convenient strategy for computing the $t$-th
component $\frak g_t$ of the desired algebra $\frak g$ as
\thetag{\ref{g}}, for a fixed integer $t=-\rho,\ldots,\varrho$.

\begin{itemize}
\item[{\bf (s1)}]  First, we construct the tangency equations
\thetag{\ref{Tangency}} of $M$ corresponding to (in general form)
holomorphic vector fields:
\[
 {\sf
X}=\sum_{j=1}^n\,Z^j(z,w)\,\partial_{z_j}+\sum_{l=1}^k\,W^l(z,w)\,\partial_{w_l}
\]
of $\frak g_t\subset\frak g$. \item[{\bf (s2)}] Now, to compute the
coefficients $Z^j(z,w)$ and $W^l(z,w)$, it is no longer necessary to
use the Taylor series \thetag{\ref{Taylor-1}} and construct and solve
the arising {\sc pde} systems as is the classical method of
\cite{Beloshapka1997,Mamai2009,MS,HAMS,Shananina2000}. Here there is
another, entirely different and much simpler way to proceed the
computation. Namely by Corollary \ref{Polynomiality}, all desired
functions $Z^j(z,w)$ and $W^l(z,w)$ are polynomials with bounded
degrees. Then, according to Proposition \ref{g(t)}(i), the desired
coefficients $Z^j(z,w)$ and $W^l(z,w)$ are weighted homogeneous
polynomials of the precise weights $[z_j]+t$ and $[w_l]+t$,
respectively. Hence, we can assume the following expressions for
them:
\begin{equation}
\label{Z,W expressions} \boxed{\aligned
Z^j(z,w)&:=\sum_{{\alpha\in\mathbb N^n\atop{\beta\in\mathbb
N^k}}\atop{[z^\alpha]+[w^\beta]=[z_j]+t}}\,{\sf
c}_{\alpha,\beta}\,.\, z^\alpha \, w^\beta,
\\
W^l(z,w)&:=\sum_{{\alpha\in\mathbb N^n\atop{\beta\in\mathbb
N^k}}\atop{[z^\alpha]+[w^\beta]=[w_l]+t}}\,{\sf
d}_{\alpha,\beta}\,.\, z^\alpha\,w^\beta, \ \ \ \ \ \ \ \
{\scriptstyle (i=1,\ldots,n, \ \ \ \ l=1,\ldots,k)},
\endaligned}
\end{equation}
for some (unknown yet) complex {\it free parameters} ${\sf
c}_{\alpha,\beta}$ and ${\sf d}_{\alpha,\beta}$. Here, by $z^\alpha$
we mean $z_1^{\alpha_1}\,z_2^{\alpha_2}\ldots z_n^{\alpha_n}$ for
$\alpha:=(\alpha_1,\ldots,\alpha_n)\in\mathbb N^n$. Furthermore, we
have:
\[
[z^\alpha]=\sum_{j=1}^n\,[z_j^{\alpha_j}]=\sum_{j=1}^n\,\alpha_j\,[z_j].
\]
Similar notations hold for $w^\beta$ and $[w^\beta]$. \item[{\bf
(s3)}] Determining the parameters ${\sf c}_{\alpha,\beta}$ and ${\sf
d}_{\alpha,\beta}$ is in fact equivalent to find the explicit
expressions of the CR-automorphisms $\sf X$ of $\frak g_t$. For this,
we should just put the already assumed expressions \thetag{\ref{Z,W
expressions}} in the $k$ tangency equations \thetag{\ref{Tangency}}
and next solve the extracted ({\it not} {\sc pde}) homogeneous linear
system of equations with the unknowns ${\sf c}_{\alpha,\beta}$ and
${\sf d}_{\alpha,\beta}$, and with the equations obtained as the
coefficients of the various powers $z^\alpha w^\beta$ of the
variables $z$ and $w$.
\end{itemize}

Let us denote by ${\sf Sys}^{t,j}$ the system of equations, mentioned
in the step ${\bf (s3)}$, associated to a CR manifold $M$ of CR
dimension and codimension $n$ and $k$, extracted from the $j$-th
tangency equation for $j=1,\ldots,k$ along the way of computing the
$t$-th component $\frak g_t$ of $\frak g=\frak{aut}_{CR}(M)$.
Furthermore, we denote by ${\sf Sys}^t$ the general system of
equations:
\[
{\sf Sys}^t:=\bigcup_{j=1}^k\,{\sf Sys}^{t,j}.
\]

\begin{Definition}
A graded Lie algebra of the form:
\[
\frak g_-:=\frak g_{-\rho}\oplus\cdots\oplus\frak g_{-1}, \ \ \ \
\rho\in\mathbb N
\]
is called {\sl fundamental} whenever it can be generated by $\frak
g_{-1}$.
\end{Definition}

In the case that the negative part $\frak g_-$ of the desired algebra
$\frak g$ is fundamental, it is even not necessary to compute the
homogeneous components $\frak g_{-\rho},\ldots,\frak g_{-2}$. In this
case, one can easily compute the $(-t)$-th components $\frak g_{-t}$
inductively as the length $t$ iterated Lie brackets $\frak
g_{-t}=[\frak g_{-1},\frak g_{-t+1}]$. There are many situations
where this occurs. For example, all Beloshapka's homogeneous models
enjoy it.

\begin{Proposition}
\label{Beloshapka-Prop} (see \cite{Beloshapka2004}, Proposition 4).
For a Beloshapka's homogeneous CR manifold $M$ as those of
\cite{Beloshapka2004}, the associated Levi-Tanaka algebra is the
negative part $\frak g_-$ of its Lie algebra of infinitesimal
CR-automorphisms. Moreover, each component $\frak g_{-m}, \,
m\in\mathbb N$, is a linear combination of brackets of degree $m$ of
vector fields in $\frak g_{-1}$.
\end{Proposition}

Then, in the case that $\frak g_-$ is fundamental, one plainly can
add the following item to the already presented strategy {\bf
(s1)-(s3)}.

\begin{itemize}
\item[{\bf (ib)}] After computing the basis elements of $\frak
g_{-1}$, to achieve each component $\frak g_{-m}$, $m=2,\ldots,\rho$
one should just compute all the length $m$ iterated brackets like:
\[
[{\sf x}_{i_1}, [{\sf x}_{i_2},[{\sf x}_{i_3},\ldots,[{\sf
x}_{i_m}]]]\ldots ], \ \ \ \ \ \ {\scriptstyle i_1<i_2< \cdots <i_m}
\]
of the basis elements ${\sf x}_\bullet$ of $\frak g_{-1}$.
\end{itemize}

\begin{Example}
\label{(1,3)-3}
 ({\it compare with Example} \ref{(1,3)}). Consider the CR-submanifold
$M\subset\mathbb C^{1+3}$ of
CR-dimension $n=1$ and codimension $k=3$ as in Example \ref{(1,3)}.
Here, we have $[z]=1, [w_1]=2$ and $[w_2]=[w_3]=3$. According to
\thetag{\ref{Tangency}}, the three fundamental tangency equations
are:
\begin{eqnarray}
\label{tangency-example}
 && 0 \equiv \big[ W^1 - \overline{W}^1 -
2i\overline{z}Z - 2iz\overline{Z}
\big]_{w=\overline{w}+\Xi(z,\overline{z},\overline{w})},
\\
\nonumber
 && 0 \equiv \big[W^2 - \overline{W}^2 - 4iz\overline{z}Z -
2i\overline{z}^2Z - 2iz^2\overline{Z} - 4iz\overline{z}\overline{Z}
\big]_{w=\overline{w}+\Xi(z,\overline{z},\overline{w})},
\\
\nonumber
 && 0 \equiv \big[W^3 - \overline{W}^3 - 4z\overline{z}Z -
2z^2\overline{Z} + 2\overline{z}^2Z + 4z\overline{z}\overline{Z}
\big]_{w=\overline{w}+\Xi(z,\overline{z},\overline{w})}.
\end{eqnarray}
First let us compute the subalgebra $\frak g_{-1}$. For this aim, we
may set the following expressions for the unknown functions
$Z(z,w_1,w_2,w_3),W^l(z,w_1,w_2,w_3)$ for $ l=1,2,3$ with their
homogeneities at their left hand sides ({\it cf.} \thetag{\ref{Z,W
expressions}}):
\[
\aligned {\tiny \boxed{ 1+(-1)}}&: \ \ \ \ Z(z,w):= {\sf p},
\\
{\tiny \boxed{ 2+(-1)}}&: \ \ \ \ W^1(z,w):= {\sf q}\,z,
\\
{\tiny \boxed{ 3+(-1)}}&: \ \ \ \ W^2(z,w):={\sf r}\,w_1+{\sf
s}\,z^2,
\\
{\tiny \boxed{ 3+(-1)}}&: \ \ \ \ W^3(z,w):={\sf t}\,w_1+{\sf
u}\,z^2,
\endaligned
\]
for some six complex functions $\sf p,q,r,s,t,u$. Putting these
expressions into the tangency equations and equating to zero the
coefficients of the appeared polynomials of $\mathbb
C[z,w_1,w_2,w_3]$, we obtain the following three systems of linear
homogeneous equations:
\[\footnotesize
\aligned {\sf Sys}^{-1,1}&=\bigg\{-2\,i\,\overline{\sf p}+{\sf q}=0,
\ \ \ \ -2\,i\,{\sf p}-\overline{\sf q}=0 \bigg\},
\\
{\sf Sys}^{-1,2}&=\bigg\{{\sf s}-2\,i\,\overline{\sf p}=0, \ \ \ \
-2\,{\sf p}-2\,\overline{\sf p}+{\sf r}=0, \ \ \ \ 2\,{\sf
p}+2\,\overline{\sf p}-{\sf r}=0\bigg\},
\\
{\sf Sys}^{-1,3}&= \bigg\{{\sf u}-2\,\overline{\sf p}=0, \ \ \ \
-2\,i\,\overline{\sf p}+2\,i\,{\sf p}+{\sf t}=0 \bigg\}.
\endaligned
\]
Solving the homogeneous linear system ${\sf Sys}^{-1}$ of all the
above equations, we can write:
\[
{\sf p}:={\sf a}+i\,{\sf b}, \ \ \ {\sf q}={\sf s}:=2({\sf b}+i\,{\sf
a}), \ \ \ {\sf r}:=4\,{\sf a}, \ \ \ {\sf t}:=4\,{\sf b}, \ \ \ {\sf
u}:=2({\sf a}-i\,{\sf b}),
\]
for two real constants $\sf a$ and $\sf b$ which brings the following
general expressions for the desired holomorphic coefficients of the
elements $\sf X\in\frak g_{-1}$ ({\it compare with} Example:
\ref{(1,3)})
\begin{equation*}
\aligned Z(z,w)&={\sf a}+i\,{\sf b},
\\
W^1(z,w)&=2\,({\sf b}+i\,{\sf a})\,z,
\\
W^2(z,w)&=2\,({\sf b}+i\,{\sf a})\,z^2+4\,{\sf a}\,w_1,
\\
W^3(z,w)&=2\,({\sf a}-i\,{\sf b})\,z^2+4\,{\sf b}\,w_1.
\endaligned
\end{equation*}
Thanks to the two free parameters $\sf a,b$ appeared in the above
expressions, we will have two infinitesimal CR-automorphisms as the
generators of $\frak g_{-1}$:
\[
\aligned {\sf
X}_1&:=\partial_z+2\,i\,z\partial_{w_1}+(2\,i\,z^2++4\,w_1)\,\partial_{w_2}+2\,z^2\,\partial_{w_3},
\\
{\sf
X}_2&:=i\,\partial_z+2\,z\,\partial_{w_1}+2\,z^2\,\partial_{w_2}+(-2\,i\,z^2+4\,w_1)\,\partial_{w_3}.
\endaligned
\]
Now, we can follow the step ${\bf (ib)}$ to seek the elements of the
homogeneous components $\frak g_{-2}$ and $\frak g_{-3}$ by computing
the iterated Lie brackets of ${\sf X}_1$ and ${\sf X}_2$ (here notice
that $M$ is one of the Beloshapka's homogeneous models and hence
$\frak g$ is fundamental. Moreover, since the maximum homogeneity of
the appearing variables is $3$ then, the minimum homogeneity in
$\frak g$ is $-3$). For $\frak g_{-2}$ we have only one generator:
\[
{\sf X_3}:=[{\sf X}_1,{\sf X}_2]=4\,\partial_{w_1}.
\]
Then, $\frak g_{-3}$ includes two basis elements
\[
\aligned {\sf X_4}&:=[{\sf X}_1,{\sf X}_3]=-4\,\partial_{w_2},
\\
{\sf X_5}&:=[{\sf X}_2,{\sf X}_3]=-4\,\partial_{w_3}.
\endaligned
\]
At present we found $5=2\times 1+3$ basis elements for the subalgebra
$\frak g_-=\frak g_{-3}\oplus\frak g_{-2}\oplus\frak g_{-1}$.
Similarly, one achieves the zeroth component $\frak g_0$ of $\frak
g$. For this, we may assume the following expressions for the
functions $Z(z,w_1,w_2,w_3),W^l(z,w_1,w_2,w_3)$ for $ l=1,2,3$ with
their homogeneities at their left hand sides:
\[
\aligned {\tiny \boxed{ 1+0}}&: \ \ \ \ Z(z,w):={\sf p_1}\,z,
\\
{\tiny \boxed{ 2+0}}&: \ \ \ \ W^1(z,w):={\sf q_1}\,w_1+{\sf
q_2}\,z^2,
\\
{\tiny \boxed{ 3+0}}&: \ \ \ \ W^2(z,w):={\sf r_1}\,w_2+{\sf
r_2}\,w_3+{\sf r_3}\,z^3+{\sf r_4}\,zw_1,
\\
{\tiny \boxed{ 3+0}}&: \ \ \ \ W^3(z,w):={\sf s_1}\,w_2+{\sf
s_2}\,w_3+{\sf s_3}\,z^3+{\sf s_4}\,zw_1.
\endaligned
\]
Again, substituting the recent expressions in the tangency equations
\thetag{\ref{tangency-example}} and equating to zero the coefficients
of the appeared polynomials of $\mathbb C[z,w_1,w_2,w_3]$, we get the
following systems of equations:
\[
\footnotesize\aligned {\sf Sys^{0,1}}&=\bigg\{{\sf q_2}=0, \ \ \ \
-{\sf p_1}-\overline{\sf p_1}+{\sf q_1}=0\bigg\},
\\
{\sf Sys^{0,2}}&=\bigg\{{\sf r_3}=0, \ \ \ \ -{\sf p_1}+\overline{\sf
p_1}+{\sf r_4}-\overline{\sf r_4}=0, \ \ \ \ \overline{\sf r_4}=0, \
\ \ \ {\sf r_1}-{\sf p_1}-2\,\overline{{\sf p_1}}+\overline{\sf
r_4}=0, \ \ \ \ {\sf r_2}=0\bigg\},
\\
{\sf Sys}^{0,3}&=\bigg\{{\sf s_3}=0, \ \ \ \ {\sf s_4}-\overline{\sf
s_4}=0, \ \ \ \ \overline{\sf s_4}=0, \ \ \ \ \overline{\sf s_4}+{\sf
s_1}+\frac{i}{2}\,{\sf p_1}-\frac{i}{2}\,\overline{\sf p_1}=0, \ \ \
\ -\overline{\sf s_4}-\overline{\sf s_1}-\frac{i}{2}\,{\sf
p_1}+\frac{i}{2}\,\overline{\sf p_1}=0,
\\
& \ \ \ \ \ \ \ \ \ \ \ \ \ \ \ \ \ \ \  -\frac{3}{2}\,{\sf
p_1}-\frac{3}{2}\,\overline{\sf p_1}+{\sf s_2}=0\bigg\}.
\endaligned
\]
Solving the linear homogeneous system ${\sf Sys}^0$ of all the above
equations, we have the solutions:
\[
\aligned {\sf p_1}&={\sf d}+i{\sf e}, \ \ \ \ {\sf q_1}=2\,{\sf d}, \
\ \ \ {\sf r_1}={\sf s_2}=3\,{\sf d}, \ \ \ \ {\sf r_2}=-{\sf e}, \ \
\ \ {\sf s_1}={\sf e},
\\
{\sf q_2}&={\sf r_3}={\sf r_4}={\sf s_3}=0,
\endaligned
\]
which implies the following expressions for the desired functions:
\begin{equation*}
\aligned Z(z,w)&=({\sf d}+i\,{\sf e})\,z,
\\
W^1(z,w)&=2\,{\sf d}\,w_1,
\\
W^2(z,w)&=3\,{\sf d}\,w_2-{\sf e}\,w_3,
\\
W^3(z,w)&={\sf e}\,w_2+3\,{\sf d}\,w_3,
\endaligned
\end{equation*}
where $\sf d$ and $\sf e$ are two real constants. Extracting the
coefficients of these two integers brings the following two tangent
vector fields, belonging to $\frak g_0$:
\[
\aligned {\sf
X}_6&:=z\,\partial_z+2\,w_1\,\partial_{w_1}+3\,w_2\,\partial_{w_2}+3\,w_3\,\partial_{w_3},
\\
{\sf X}_7&:=i\,z\,\partial_z-w_3\,\partial_{w_2}+w_2\,\partial_{w_3}.
\endaligned
\]
\end{Example}

\begin{Remark}
It is worth noting that for computing each subspace $\frak g^{(t)}$,
although one can achieve its basis elements by computing the
corresponding components $\frak g_s,\, t=-\rho,\ldots,t$, it is also
possible to adopt the above strategy {\bf (s1)}-{\bf (s3)} by
modifying the assumed expressions of the functions $Z^j(z,w)$ and
$W^l(z,w)$ as follows ({\it cf.} Proposition \ref{g(t)}$(ii)$):
\begin{equation}
\label{Z,W expressions for g(t)} \aligned
Z^j(z,w)&:=\sum_{{\alpha\in\mathbb N^n\atop{\beta\in\mathbb
N^k}}\atop{[z^\alpha]+[w^\beta]\leqslant [z_j]+t}}\,{\sf
c}_{\alpha,\beta}\,.\, z^\alpha \, w^\beta,
\\
W^l(z,w)&:=\sum_{{\alpha\in\mathbb N^n\atop{\beta\in\mathbb
N^k}}\atop{[z^\alpha]+[w^\beta]\leqslant [w_l]+t}}\,{\sf
d}_{\alpha,\beta}\,.\, z^\alpha\,w^\beta, \ \ \ \ \ \ \ \
{\scriptstyle (i=1,\ldots,n, \ \ \ \ l=1,\ldots,k)}.
\endaligned
\end{equation}
In particular when $M$ has rigidity, we can obtain the sought algebra
$\frak g$ by setting $t=0$.
\end{Remark}

\subsection{Finding the maximum homogeneity $\varrho$}

For an arbitrary homogeneous CR manifold, represented as
\thetag{\ref{model}}, so far we have provided an effective way to
compute homogeneous components $\frak g_t$ of the graded desired
algebra $\frak g:=\frak{aut}_{CR}(M)$ of the form:
\[
\frak g:=\frak g_{-\rho}\oplus\cdots\oplus\frak
g_0\oplus\cdots\oplus\frak g_{\varrho}.
\]
We also find that the value of $\rho$ in this gradation is equal to
the maximum weight $[w_k]$ appearing among the complex variables. The
only not-yet-fixed problem is to answer how much we have to compute
the homogeneous components $\frak g_t$ to arrive at the last one
$\frak g_\varrho$. Here, we do not aim to find the precise value of
$\varrho$ but\,\,---\,\,in an algorithmic point of view\,\,---\,\,it
suffices to find a criterion to stop the computations. For this aim,
we employ the transitivity of the Lie algebra $\frak g$. For every
homogeneous CR manifold $M$, the associated Lie algebra of its
infinitesimal CR-automorphisms is transitive:

\begin{Definition}
A graded Lie algebra $\frak g$ as above is called {\sl transitive}
whenever for each element $\sf x\in\frak g_t$ with $t\geqslant  0$, the
equality $[{\sf x},\frak g_-]=0$ implies that ${\sf x}=0$. In the
case that $\frak g_-$ is fundamental then, the transitivity means
that for any $\sf x$ as above, the equality $[{\sf x},\frak
g_{-1}]=0$ implies that ${\sf x}=0$.
\end{Definition}

\begin{Proposition}
Consider a transitive graded algebra $\frak g$ as above. For each
integer $t\geqslant 0$, if $\frak g_t=\frak g_{t+1}=\cdots=\frak
g_{t+\rho-1}\equiv 0$ then we have $\frak g_{t+\rho}=0$. Moreover, if
$\frak g$ is also fundamental then the equality $\frak g_t=0$ implies
independently that $\frak g_{t+1}=0$.
\end{Proposition}

\proof Assume the following gradation for the transitive algebra
$\frak g$:
\[
\frak g:=\frak g_-\oplus\frak g_0\oplus\frak
g_1\oplus\ldots\oplus\frak g_{t-1}\oplus \underbrace{0}_{\frak g_t}
\oplus\underbrace{0}_{\frak
g_{t+1}}\oplus\ldots\oplus\underbrace{0}_{\frak
g_{t+\rho-1}}\oplus\frak g_{t+\rho}\oplus\ldots
\]
and let ${\sf x}\in\frak g_{t+\rho}$. According to the inequality
$[\frak g_i,\frak g_j]\subset\frak g_{i+j}$, we have:
\[
\aligned &[{\sf x},\frak g_{-1}]\subset\frak g_{t+\rho-1}=0,
\\
&[{\sf x},\frak g_{-2}]\subset\frak g_{t+\rho-2}=0,
\\
&\ \ \ \ \ \ \ \ \  \ \ \vdots
\\
&[{\sf x},\frak g_{-\rho}]\subset\frak g_{t+\rho-\rho}=0,
\endaligned
\]
which implies that $[{\sf x},\frak g_-]=0$. Now, the transitivity of
$\frak g$ immediately implies that ${\sf x}=0$. For the second part
of the assertion, similarly assume the following gradation for $\frak
g$:
\[
\frak g:=\frak g_-\oplus\frak g_0\oplus\frak
g_1\oplus\ldots\oplus\frak g_{t-1}\oplus \underbrace{0}_{\frak g_t}
\oplus\frak g_{t+1}\oplus\ldots,
\]
and let $\sf x\in\frak g_{t+1}$. Consequently we have:
\[
[{\sf x},\frak g_{-1}]\subset\frak g_{t+1-1}=\frak g_{t}=0.
\]
Again, the definition of fundamental transitive algebras immediately
implies that ${\sf x}=0$. This completes the proof.
\endproof

Accordingly, for a homogeneous CR manifold $M$ and to realize how
much we have to compute the homogeneous components of $\frak g$ to
arrive at the last weighted homogeneous component $\frak g_\varrho$
we can apply the following plain strategy:
\begin{itemize}
\item[$\bullet$] {\bf When $\frak g_-$ is fundamental.} Compute the
homogeneous components $\frak g_t$ of $\frak g$ as much as it appears
the first trivial component.
 \item[$\bullet$] {\bf When $\frak g_-$
is not fundamental.} Compute the homogeneous components $\frak g_t$
of $\frak g$ as much as they appear $\rho$ successive trivial
components.
\end{itemize}

\begin{Example}
In Example \ref{(1,3)-3}, we computed the negative part $\frak
g_-=\frak g_{-3}\oplus\frak g_{-2}\oplus\frak g_{-1}$ and also zeroth
component $\frak g_0$ of the graded algebra $\frak
g:=\frak{aut}_{CR}(M)$ associated to the presented homogeneous
CR manifold $M\subset\mathbb C^{1+3}$. Here, let us finalize
computation of the desired algebra. Proceeding further in this
direction, now let us compute the next component $\frak g_1$. In this
case, we can set the following expressions for four desired
coefficients $Z(z,w)$ and $W^l(z,w)\, l=1,2,3$, with their weighted
homogeneities at their left hand sides ({\it cf.} \thetag{\ref{Z,W
expressions}})
\[
\aligned {\tiny \boxed{ 1+1}}&: \ \ \ \ Z(z,w):= {\sf a}_1\,w_1+{\sf
a}_2\,z^2,
\\
{\tiny \boxed{ 2+1}}&: \ \ \ \ W^1(z,w):= {\sf a}_3\,w_2+{\sf
a}_4\,w_3+{\sf a}_5\,z^3+{\sf a}_6\,zw_1,
\\
{\tiny \boxed{ 3+1}}&: \ \ \ \ W^2(z,w):={\sf a}_7\,w_1^2+{\sf
a}_8\,z^2w_1+{\sf a}_9\,z^4+{\sf a}_{10}\,zw_2+{\sf a}_{11}\,zw_3,
\\
{\tiny \boxed{ 3+1}}&: \ \ \ \ W^3(z,w):={\sf a}_{12}\,w_1^2+{\sf
a}_{13}\,z^2w_1+{\sf a}_{14}\,z^4+{\sf a}_{15}\,zw_2+{\sf
a}_{16}\,zw_3.
\endaligned
\]
Now, we should check these predefined expressions in the tangency
equations \thetag{\ref{tangency-example}}. This gives us the total
system ${\sf Sys}^1=\bigcup_{j=1}^3 {\sf Sys}^{1,j}$ as follows:
\begin{equation*}
\aligned {\sf Sys}&^1=
\\
&\left\{
\begin{array}{l}
{\sf a}_5=0, \ \ \ i\,\big({\sf a}_6+{\sf a}_3-{\sf a}_2=0\big)+{\sf
a}_4=0, \ \ \ i\,\big({\sf a}_3-\overline{\sf a}_2\big)-{\sf
a}_4+2\,{\sf a}_1=0, \ \ \ -2i\,\overline{\sf a}_1+{\sf a}_6=0,
\\
\overline{\sf a}_5=0, \ \ \ -2i\,{\sf a}_1-\overline{\sf a}_6=0, \ \
\ {\sf a}_3-\overline{\sf a}_3=0, \ \ \ {\sf a}_4-\overline{\sf
a}_4=0, \ \ \ {\sf a}_9=0, \ \ \ i\,\big({\sf a}_{10}-2\,{\sf
a}_2+{\sf a}_8\big)+{\sf a}_{11}=0,
\\
i\,\big({\sf a}_{10}-{\sf a}_2-\overline{\sf a}_2\big)+8\,{\sf
a}_1-2\,{\sf a}_{11}-4\,{\sf a}_7=0, \ \ \ -2i\,\overline{\sf

a}_1+{\sf a}_8=0, \ \ \ {\sf a}_1-i\,\overline{\sf a}_2=0, \ \ \
{\sf a}_{10}=0,
\\
{\sf a}_{11}=0, \ \ \ i\,\big({\sf a}_7-{\sf a}_1-\overline{\sf
a}_1=0\big)=0, \ \ \ \overline{\sf a}_9=0, \ \ \ -\overline{\sf

a}_8-2i\,{\sf a}_1=0, \ \ \ \overline{\sf a}_{10}=0, \ \ \
\overline{\sf a}_{11}=0,\ \ \  {\sf a}_{14}=0,
\\
{\sf a}_7-\overline{\sf a}_7=0, \ \ \ i\,\big({\sf a}_{13}+{\sf
a}_{15}\big)+{\sf a}_{16}-2\,{\sf a}_2=0, \ \ \ -2\,{\sf a}_{12}+{\sf
a}_2-{\sf a}_{16}-\overline{\sf a}_2+i\,\big({\sf a}_{15}-4\,{\sf
a}_1\big)=0,
\\
{\sf a}_{13}-2\,\overline{\sf a}_1=0, \ \ \ \overline{\sf
a}_2+i\,{\sf a}_1=0, \ \ \ i\,{\sf a}_{12}-{\sf a}_1+\overline{\sf
a}_1=0, \ \ \ {\sf a}_{15}=0, \ \ \ {\sf a}_{16}=0, \ \ \
\overline{\sf a}_{14}=0,
\\
-\overline{\sf a}_{13}+2\,{\sf a}_1=0, \ \ \ -\overline{\sf

a}_{15}=0, \ \ \ -\overline{\sf a}_{16}=0, \ \ \ {\sf
a}_{12}-\overline{\sf a}_{12}=0
\end{array}
\right\}.
\endaligned
\end{equation*}
This system has only the trivial solution ${\sf a}_1=\cdots={\sf
a}_{16}=0$ which means that we have $\frak g_1=0$. Furthermore as we
know, the desired algebra $\frak g$ is fundamental which guarantees
that the next components are trivial, too. Summing up the results of
this example with those of Example \ref{(1,3)-3}, one finds the
sought 7-dimensional Lie algebra of infinitesimal CR-automorphisms
associated to $M$ as the gradation:
\[
\frak g=\frak g_{-3}\oplus\frak g_{-2}\oplus\frak g_{-1}\oplus\frak
g_0,
\]
with $\frak g_{-3}=\langle{\sf X}_4,{\sf X}_5\rangle$, with $\frak
g_{-2}=\langle{\sf X}_3\rangle$, with $\frak g_{-1}=\langle{\sf
X}_1,{\sf X}_2\rangle$ and with $\frak g_{0}=\langle{\sf X}_6,{\sf
X}_7\rangle$ with the Lie commutators displayed in the following
table:
\medskip
\begin{center}
\begin{tabular} [t] { c | c c c c c c c }
& ${\sf X}_5$ & ${\sf X}_4$ & ${\sf X}_3$ & ${\sf X}_2$ & ${\sf X}_1$
& ${\sf X}_6$ & ${\sf X}_7$
\\
\hline ${\sf X}_5$ & $0$ & $0$ & $0$ & $0$ & $0$ & $3{\sf X}_5$ &
$-{\sf X}_4$
\\
${\sf X}_4$ & $*$ & $0$ & $0$ & $0$ & $0$ & $3{\sf X}_4$ & ${\sf
X}_5$
\\
${\sf X}_3$ & $*$ & $*$ & $0$ & $4{\sf X}_5$ & $4{\sf X}_4$ & $2{\sf
X}_3$ & $0$
\\
${\sf X}_2$ & $*$ & $*$ & $*$ & $0$ & $-4{\sf X}_3$ & ${\sf X}_2$ &
$-{\sf X}_1$
\\
${\sf X}_1$ & $*$ & $*$ & $*$ & $*$ & $0$ & ${\sf X}_1$ & ${\sf X}_2$
\\
${\sf X}_6$ & $*$ & $*$ & $*$ & $*$ & $*$ & $0$ & $0$
\\
${\sf X}_7$ & $*$ & $*$ & $*$ & $*$ & $*$ & $*$ & $0$ %
\end{tabular}
\end{center}
One observes that the achieved algebra is exactly that of Example
\ref{(1,3)}.
\end{Example}

One can finds another computation of the Lie algebra achieved in the
above example via the classical method of solving the arisen {\sc
pde} system in \cite{5-cubic}. Comparing the above process with that
of this paper clarify the effectiveness of the prepared algorithm.

\section{Summing up the results}

\label{section main strategy}

Here, let us gather the results obtained so far
to provide an algorithm for computing the sought Lie algebras of
infinitesimal CR-automorphisms associated to the holomorphically
nondegenerate homogeneous CR manifolds, represented as
\thetag{\ref{model-model}}. The strategy introduced in subsection
\ref{strategy-1} enabled one to compute separately the homogeneous
components $\frak g_t, \, t=-\rho,\ldots,\varrho$ of the graded
algebra $\frak g$ of infinitesimal CR-automorphisms of such
CR manifolds as:
\begin{equation}
\label{gg} \frak g:=\frak g_{-\rho}\oplus\cdots\oplus\frak
g_{-1}\oplus\frak g_0\oplus\cdots\oplus\frak g_\varrho.
\end{equation}
One may follow the following two points for computing the desired
algebras associated to the under consideration homogeneous
CR manifolds:

\begin{itemize}
\item[\bf{Point 1}] Executing three steps $\bf{(s1)-(s2)-(s3)}$
introduced in subsection \ref{strategy-1} and finding homogeneous
components $\frak g_t$ successively. In particular if $\frak g_-$ is
fundamental then one can execute the step {\bf (ib)}. \item[\bf{Point
2}] In the above gradation, the minimum homogeneity $-\rho$ is equal
to $[w_k]$ where $w_k$ has the maximum homogeneity among the complex
variables appearing in \thetag{\ref{model-model}}. Moreover, it
suffices to compute successively the homogeneous components $\frak
g_{0},\frak g_{1},\ldots$ as much as we find $\rho$ successive
trivial algebras. In particular if $\frak g_-$ is fundamental, we can
terminate the computations as much as we find first trivial
component.
\end{itemize}

\begin{Remark}
In \cite{Beloshapka2004}, Beloshapka called his introduced
CR manifolds by {\it nice} universal CR-models. These manifolds have
the ability of enjoying all properties, required for launching the
presented algorithm. They are of finite type, holomorphically
nondegenerate, generic and real analytic with the graded Lie algebras
of infinitesimal CR-automorphisms their negative parts are
fundamental. Having such properties, we shall confirm that these
models are really deserved to be called by the phrase {\it nice}.
\end{Remark}

Let us conclude this paper by computing Lie algebras of infinitesimal
CR-automorphisms associated to the following CR-manifold:

\begin{equation}
\label{New-Model} \aligned \mathcal M:=\left\{
\begin{array}{l}
w_1-\overline w_1=2i\,z\overline z,
\\
w_2-\overline w_2=2i\,z\overline z\,\big(z+\overline z\big),
\\
w_3-\overline w_3=2i\,z\overline
z\,\big(z^2+\textstyle{\frac{3}{2}}\,z\overline z+\overline z^2\big).
\end{array}
\right.
\endaligned
\end{equation}
It is worth noting that this CR-manifold\,\,---\,\,which admits some
interesting features \cite{MPS, 5-cubic}\,\,---\,\,does not belong to the
class of CR-models introduced by Beloshapka \cite{Beloshapka2004}.
Before computations, let us show that $\mathcal M$ is homogeneous.

\begin{Proposition}
The already introduced CR-manifold $\mathcal M\subset\mathbb C^{1+3}$
is homogeneous.
\end{Proposition}

\proof To prove the assertion, we show that for each arbitrary point
$(p,q_1,q_2,q_3)$ of $\mathcal M$, there is a CR-automorphism mapping
the origin to it. For this aim, it suffices to prove that after
applying the four replacements $z\mapsto z+p$ and $w_j\mapsto
w_j+q_j, \, j=1,2,3$, one can transform the model into its initial
form by means of some certain changes of coordinates. At first, let
us consider the effect of these transformations on the third
equation. In fact, it gives:
\begin{equation}
\label{after-transform}\aligned w_3-\overline w_3&+q_3-\overline
q_3=2i\,z\overline z\,\big(z^2-3\,z\overline z+\overline z^2\big)+
\\
&+3i\,p^2\overline p^2+2i\,p^3\overline p+2i\,p\overline
p^3+(6i\,\overline p^2p+2i\,\overline p^3+6i\,\overline
pp^2)z+(6i\,\overline pp^2+6i\,\overline p^2p+ 2i\,p^3)\overline z+
\\
&+(6i\,\overline p^2+6i\,p^2+12i\,p\overline p)z\overline z+
(3i\,\overline p^2+6i\,p\overline p)z^2+(6i\,p\overline
p+3i\,p^2)\overline z^2+
\\
&+(6i\,p+6i\,\overline p)z^2\overline z+(6i\,\overline
p+6i\,p)z\overline z^2+2i\,\overline p z^3+2i\,p\overline z^3.
\endaligned
\end{equation}
It is possible to kill all pluriharmonic terms by the change of
variables $w_3\mapsto w_3+P_3(z)$ with:
\[
\aligned P_3(z)=-q_3+\textstyle{\frac{3}{2}}i\,p^2\overline
p^2+2i\,p^3\overline p+\big(6i\,p\overline p^2+2i\,\overline
p^3+6i\,p^2\overline p\big)\,z+\big(3i\,\overline p^2+6i\,p\overline
p\big)\,z^2+2i\overline p\,z^3.
\endaligned
\]
 After such change of
variables we will have:
\[\aligned
w_3-\overline w_3&=2i\,z\overline z\,\big(z^2-3\,z\overline
z+\overline z^2\big)
\\
&+(6i\overline p^2+6ip^2+12ip\overline p)z\overline z+
\\
&+(6ip+6i\overline p)z^2\overline z+(6i\overline p+6ip)z\overline
z^2.
\endaligned
\]
Replacing the term $z\overline z$ by the equal expression
$\frac{1}{2i}\,(w_1-\overline w_1)$ at the second line, one plainly
verifies that it is also possible to eliminate this line by the
change of coordinates $w_3\mapsto w_3+\big(3\,\overline
p^2+3\,p^2+6\,p\overline p\big)\,w_1$ and obtain:
\[\aligned
w_3-\overline w_3&=2i\,z\overline z\,\big(z^2-3\,z\overline
z+\overline z^2\big)
\\
&+(6ip+6i\overline p)z^2\overline z+(6i\overline p+6ip)z\overline
z^2.
\endaligned
\]
Now it remains to eliminate the second line of the above expression
again by some holomorphic change of coordinates. First, one should
notice that according to the defining equations of the model, the
second line can be represented into the form:
\[
(6ip+6i\overline p)z^2\overline z+(6i\overline p+6ip)z\overline
z^2=(3\overline p+3p)\,\big(w_2-\overline w_2\big).
\]
Hence, to eliminate this line from the last expression, it suffices
to use the holomorphic change of coordinates $w_3\mapsto
w_3+(3\overline p+3p)\,w_2$. This convert the third expressions into
the initial form:
\[\aligned
w_3-\overline w_3&=2i\,z\overline z\,\big(z^2-3\,z\overline
z+\overline z^2\big),
\endaligned
\]
as desired. Simpler procedure works in the cases of two first
defining equations.
\endproof

\begin{Example}
Now, let us compute the Lie algebra $\frak
g:=\frak{aut}_{CR}(\mathcal M)$ associated to the CR-manifold
$\mathcal M$ defined as \thetag{\ref{New-Model}}. First, one notices
that we have the following weights of the appearing complex
variables:
\[
[z]=1, \ \ \ [w_1]=2, \ \ \ [w_2]=3, \ \ \ [w_3]=4.
\]
Hence the minimum homogeneity of the homogeneous components will be
$-\rho=-4$. Here, an infinitesimal CR-automorphism is of the form:
\[
{\sf
X}:=Z(z,w)\,\partial_z+W^1(z,w)\,\partial_{w_1}+W^2(z,w)\,\partial_{w_2}+W^3(z,w)\,\partial_{w_3},
\]
enjoying the following three fundamental tangency equations:
\begin{equation}\label{tangency-2}
\aligned 0&\equiv\big[W^1-\overline W^1-2i\,\overline z
Z-2i\,z\overline Z\big]_{\mathcal M},
\\
0&\equiv\big[W^2-\overline W^2-4i\,z\overline z Z-2i\,\overline z^2
Z-2i\,z^2\overline Z-4i\,z\overline z\overline Z\big]_{\mathcal M},
\\
0&\equiv\big[W^3-\overline W^3-6i\,z^2\overline z Z-6i\,z\overline
z^2 Z-2i\,\overline z^3 Z-2i\,z^3\overline Z-6i\,z^2\overline
z\overline Z-6i\,z\overline z^2 \overline Z\big]_{\mathcal M}.
\endaligned
\end{equation}
Let us start by computing the negative part $\frak g_-$. For the
$(-1)$-th component $\frak g_{-1}$, the sought coefficients are of
the forms:
\begin{equation*}\aligned
\left[
\begin{array}{l}
Z(z,w):={\sf a}_1,
\\
W^1(z,w):={\sf a}_2\,z,
\\
W^3(z,w):={\sf a}_3\,z^2+{\sf a}_4\,w_1,
\\
W^4(z,w):={\sf a}_5\,z^3+{\sf a}_6\,zw_1+{\sf a}_7\,w_2.
\end{array}
\right.
\endaligned
\end{equation*}
Checking these predefined polynomials into the tangency equations
\thetag{\ref{tangency-2}} gives the following system:
\begin{equation*}
\aligned {\sf Sys}^{-1}:=\left\{
\begin{array}{l}
{\sf a}_2-2i\,\overline{\sf a}_1=0 \ \ \ -4i\,{\sf

a}_1-4i\,\overline{\sf a}_1+2i\,\overline {\sf a}_4=0 \ \ \ {\sf
a}_3-2i\,\overline{\sf a}_1=0
\\
{\sf a}_4-\overline{\sf a}_4=0 \ \ \ {\sf a}_5-2i\,\overline{\sf
a}_1=0 \ \ \ {\sf a}_6=0 \ \ \ {\sf a}_7-\overline{\sf a}_7=0
\end{array}
\right\}
\endaligned
\end{equation*}
which has the solution:
\[\aligned
{\sf a}_1&:={\sf a}+i\,{\sf b}, \ \ \ {\sf a}_2=2\,{\sf b}+2i\,{\sf
a}, \ \ \ {\sf a}_3=2\,{\sf b}+2i\,{\sf a}, \ \ \ {\sf a}_4=4\,{\sf
a},
\\
{\sf a}_5&=2\,{\sf b}+2i\,{\sf a}, \ \ \ {\sf a}_6=0, \ \ \ {\sf
a}_7=6\,{\sf a}, \ \ \ \ \ \ \ \ \ \ \  \ {\scriptstyle ({\sf
a,b}\,\in\,\mathbb R)}.
\endaligned
\]
Consequently, the sought homogeneous component $\frak g_{-1}$ is
2-dimensional with the generators:
\[\aligned
{\sf
X}_1&=\partial_z+2iz\,\partial_{w_1}+2iz^2\,\partial_{w_2}+4w_1\,\partial_{w_2}+2iz^3\,\partial_{w_3}+6w_2\,\partial_{w_3},
\\
{\sf
X}_2&=i\,\partial_z+2z\,\partial_{w_1}+2z^2\,\partial_{w_2}+2z^3\,\partial_{w_3}.
\endaligned
\]
Similar (and even simpler) computations give three vector fields:
\[\aligned
{\sf X}_3&:=\partial_{w_1},
\\
{\sf X}_4&:=\partial_{w_2},
\\
{\sf X}_5&:=\partial_{w_3},
\endaligned
\]
of homogeneities $-2,-3,-4$, respectively. So far, we have computed
the negative component:
\[
\frak g_-:=\frak g_{-4}\oplus\frak g_{-3}\oplus\frak
g_{-2}\oplus\frak g_{-1}
\]
with $\frak g_{-1}:=\langle {\sf X}_1,{\sf X}_2\rangle$, with $\frak
g_{-2}:=\langle {\sf X}_3\rangle$, with $\frak g_{-3}:=\langle {\sf
X}_4\rangle$ and with $\frak g_{-4}:=\langle {\sf X}_5\rangle$. One
can see all the possible Lie commutators of these generators in the
table presented below. From this table, one easily verifies that the
negative part $\frak g_-$ of $\frak g$ is in fact {\it fundamental}.
Hence, we have to compute nonnegative components as much as we
encounter first trivial one. Now, we have to continue with computing
$\frak g_0$. In this case, the sought coefficients are of the form:
\begin{eqnarray*}
\left[
\begin{array}{l}
Z := {\sf a}_1\, z\\
W_1 := {\sf a}_2\, w_1+{\sf a}_3\, z^2\\
 W_2 := {\sf a}_4\, w_2+{\sf a}_5 \,z^3+{\sf a}_6 \,z w_1\\
 W_3
:= {\sf a}_7\, w_3+{\sf a}_8\, w_1^2+{\sf a}_9\, z^2 w_1+{\sf
a}_{10}\, z^4+{\sf a}_{11}\, z w_2.
\end{array}
\right.
\end{eqnarray*}
Checking these predefined functions in the tangency equations
\thetag{\ref{tangency-example}} gives the following complex system:
\begin{equation*}
\aligned {\sf Sys}^0&=
\\
& \left\{
\begin{array}{l}
{\sf a}_3=0, \ \ \  -2i\, {\sf a}_1-2i\, \overline{{\sf a}}_1+2i\,
{\sf a}_2=0, \ \ \ -\overline{{\sf a}}_3=0, \ \ \ -\overline{{\sf

a}}_2+{\sf a}_2=0, \ \ \
 {\sf a}_5=0,
 \\
-4i\, {\sf a}_1-2i\, \overline{{\sf a}}_1+2i\,{\sf a}_4+2i\, {\sf
a}_6=0, \ \ \ -4i\,\overline{{\sf a}}_1-2i\, {\sf a}_1+2i\, {\sf
a}_4=0, \ \ \  {\sf a}_6-\overline{{\sf a}}_5=0, \ \ \
 -\overline{{\sf a}}_6=0,
 \\
    {\sf a}_{10}-2i\, \overline{{\sf a}}_1+2i\, {\sf a}_9+2i\, {\sf a}_7+2i\,{\sf
a}_{11}-6i\,{\sf a}_1=0,
  \ \ \  3i\,{\sf a}_7+2i\, {\sf a}_{11}-4{\sf a}_8-6i\, {\sf a}_1-6i\,
\overline{{\sf a}}_1=0,
  \\
{\sf a}_4-\overline{{\sf a}}_4=0, \ \ \   {\sf a}_9=0, \ \ \
2i\,{\sf a}_7-6i\,\overline{{\sf a}}_1-2i\,{\sf a}_1=0, \ \ \
4i\,{\sf a}_8=0, \ \ \  {\sf a}_{11}=0, \ \ \ -\overline{{\sf

a}}_{10}=0, \ \ \
\\
-\overline{{\sf a}}_9=0, \ \ \  -\overline{{\sf a}}_{11}=0, \ \ \
{\sf a}_8-\overline{{\sf a}}_8=0, \ \ \ {\sf a}_7-\overline{{\sf
a}}_7=0
\end{array}
\right\}.
\endaligned
\end{equation*}
This system has the solution:
\[
{\sf a}_1={\sf a}, \ \ \ {\sf a}_2=2\,{\sf a}, \ \ \ {\sf
a}_4=3\,{\sf a}, \ \ \ {\sf a}_7=4\,{\sf a}, \ \ \ {\sf a}_3={\sf
a}_5={\sf a}_6={\sf a}_8={\sf a}_9={\sf a}_{10}={\sf a}_{11}\equiv 0
\]
for some real number ${\sf a}$. Therefore, $\frak g_0$ is
1-dimensional with the generator:
\[
{\sf
X}_0=z\,\partial_z+2w_1\,\partial_{w_1}+3w_2\,\partial_{w_2}+4w_3\,\partial_{w_3}.
\]
This nonnegative component was not trivial; hence we have to proceed
by computing the next component $\frak g_1$. Similar computations
that we do not present them for saving space shows that this
component {\it is trivial}. Then, according to the fundamentality of
$\frak g_-$ we can terminate the computations. Consequently, the
sought graded algebra $\frak g$ is of the form:
\[
\frak g:=\frak g_{-4}\oplus\frak g_{-3}\oplus\frak g_{-2}\oplus\frak
g_{-1}\oplus\frak g_0
\]
with the negative components as above, with $\frak g_0=\langle{\sf
X}_0\rangle$ and with the table of commutators displayed as follows:

\medskip
\begin{center}
\begin{tabular} [t] { c | c c c c c c }
& ${\sf X}_0$ & ${\sf X}_1$ & ${\sf X}_2$ & ${\sf X}_3$ & ${\sf X}_4$
& ${\sf X}_5$
\\
\hline ${\sf X}_0$ & $0$ & $-{\sf X}_1$ & $-{\sf X}_2$ & $-2{\sf
X}_3$ & $-3{\sf X}_4$ & $-4{\sf X}_5$
\\
${\sf X}_1$ & $*$ & $0$ & $-4{\sf X}_3$ & $-4{\sf X}_4$ & $6{\sf
X}_5$ & $0$
\\
${\sf X}_2$ & $*$ & $*$ & $0$ & $0$ & $0$ & $0$
\\
${\sf X}_3$ & $*$ & $*$ & $*$ & $0$ & $0$ & $0$
\\
${\sf X}_4$ & $*$ & $*$ & $*$ & $*$ & $0$ & $0$
\\
${\sf X}_5$ & $*$ & $*$ & $*$ & $*$ & $*$ & $0$
\end{tabular}
\end{center}
\end{Example}

\section{Parametric defining equations
\\
and Gr\"obner systems}
\label{Comprehensive}

Actually, one of the main\,\,---\,\,somehow hidden\,\,---\,\,obstacles appearing among the computations
arises when the set of defining equations includes some certain {\it parametric} polynomials.
This case is quite usual as one observes in \cite{Beloshapka2004,Beloshapka1997,Mamai2009,Shananina2000}. To treat such cases, we suggest the modern and effective concept of {\sl comprehensive Gr\"obner systems} \cite{Weis, kapur,kapur2,monts1,monts3} which enables us to consider and solve (linear) parametric systems appearing among the computations.

To begin, let $\mathbb K$ be a field and let ${\bf a} := a_1,\ldots,a_t$ and ${\bf x}:=x_1,\ldots,x_n$ be two certain sequences of parameters and variables, respectively. Naturally, we call the ring:
 \[
 \mathbb K[{\bf a}][{\bf x}]:=\left\{\sum_{i=1}^{m} p_{\alpha_i} x_1^{\alpha_{i1}}\cdots x_n^{\alpha_{in}}\,|\, p_{\alpha_i}\in \mathbb K[{\bf a}], \, \alpha_{ij}\in\mathbb N\cup\{0\}\right\}
 \]
 the {\sl parametric polynomial ring} over $\mathbb K$ with parameters ${\bf a}$ and variables ${\bf x}$. Let $P$ be a set of parametric polynomials which generates the parametric  ideal ${\mathcal I}$. Obviously, the solutions of the {\sl parametric system} $P=0$ depend on the extant parameters. The main idea behind the modern concept of {\sl comprehensive Gr\"obner systems} is to treat such solutions by discussing the values of the parameters of the system defined by $\mathcal I$. This concept provides some effective and powerful tools in which enables one to divide the space of parameters into a finite number of partitions for which the {\it general form} of solutions arising from each partition is unique.

\begin{Definition}
\label{defGS}
Let ${\mathcal I}\subset \mathbb K[{\bf a}][{\bf x}]$ be a parametric ideal, $\overline{\mathbb K}$ be the algebraic closure of $\mathbb K$ and $\prec$ be a monomial ordering on ${\bf x}$. Then the set:
\[
{\bf G}({\mathcal I})=\{(E_i,N_i,G_i) \mid i=1,\ldots,\ell\} \subset \mathbb K[{\bf a}]\times \mathbb K[{\bf a}]\times \mathbb K[{\bf a}][{\bf x}]
\]
is called a comprehensive  Gr\"obner system for $\mathcal I$ if for each homomorphism $\sigma_{(\lambda_1,\ldots,\lambda_t)}  : \mathbb K[{\bf a}][{\bf x}] \longrightarrow   \overline{\mathbb K}[{\bf x}]$, associated to a $t$-tuple $(\lambda_1,\ldots,\lambda_t)\in \overline{\mathbb K}^t$ and defined by:
\begin{eqnarray*}
\begin{array}{cccc}
\sum_{i=1}^{m} p_{\alpha_i}(a_1,\cdots,a_t)\, x_1^{\alpha_{i1}}\cdots x_n^{\alpha_{in}} &\mapsto & \sum_{i=1}^{m} p_{\alpha_i}(\lambda_1,\ldots,\lambda_t)\, x_1^{\alpha_{i1}}\cdots x_n^{\alpha_{in}},\\
\end{array}
\end{eqnarray*}
 there exists a pair $(E_i,N_i)$ with $(\lambda_1,\ldots,\lambda_t) \in { V}(E_i) \setminus { V}(N_i)$ such that $\sigma(G_i)$ is a Gr\"obner basis for $\sigma(I)$ with respect to $\prec$. Here by $V(E_i)$ and $V(N_i)$ we mean the algebraic varieties associated to the polynomial sets $E_i$ and $N_i$. In this case, $E_i$ and $N_i$ are called {\sl null} and {\sl non-null conditions}, respectively.
\end{Definition}

Remark that, by \cite[Theorem 2.7]{Weis}, every parametric ideal possesses a ({\em finite}) comprehensive Gr\"obner system, however, by Definition \ref{defGS}, we can observe that such a system may be not unique.
The concept of Gr\"obner systems was introduced first by Weispfenning in 1992 \cite{Weis}. Later on, Montes \cite{monts1} proposed {\sc DisPGB} algorithm for computing Gr\"obner systems. In 2006, Sato and Suzuki \cite{susa}  provided an important improvement for computing Gr\"obner systems by doing only computation of the reduced Gr\"obner bases in polynomial rings over ground fields. Furthermore, Montes and Wibmer in \cite{monts3}  presented the {\sc Gr\"obnerCover} algorithm which computes a finite partition of the parameter space into locally closed subsets together with polynomial data  from which the {\it reduced} Gr\"obner basis for a given values of parameters can immediately be determined. Kapur, Sun and Wang \cite{kapur,kapur2} in 2010 and 2013 suggested two new algorithms for computing Gr\"obner systems  by combining Weispfenning's algorithm with Suzuki and Sato's.

It is worth noting that if ${ V}(E_i) \setminus { V}(N_i) = \emptyset$, for some $i$, then the triple $(E_i,N_i,G_i)$ is useless and it must be omitted from the comprehensive Gr\"obner system. In this case, the pair $(E_i,N_i)$ is called {\sl inconsistent}. It is known that inconsistency occurs if and only if $N_i \subset \sqrt{\langle E_i \rangle}$ and thus we need to an efficient radical membership test to determine it.

In the recently published paper \cite{kapur2}, Kapur, Sun and Wang introduced an effective algorithm to compute comprehensive Gr\"obner system of a parametric polynomial ideal. This algorithm which is called by {\sc PGB} uses a new and efficient radical membership criterion based on linear algebra methods. To the best of our knowledge, it is the most powerful algorithm of computing comprehensive Gr\"obner systems introduced so far and it is for this reason that we prefer to employ this algorithm in our computations.

Besides the deep theory encompassing this subject, the concept of comprehensive Gr\"obner bases provides some effective tools to consider and to solve parametric systems by decomposing the space of the extant parameters. To illustrate this ability let us borrow the following example from \cite{kapur2}.

\begin{Example}
Consider the following parametric polynomial system in $\mathbb{C}[a,b,c][x,y]$:
\begin{eqnarray*}
\Sigma:\left\{
\begin{array}{lll}
ax-b&=&0\\
by-a&=&0\\
cx^2-y&=&0\\
cy^2-x&=&0.
\end{array}
\right.
\end{eqnarray*}
Choosing the graded reverse lexicographical ordering $y\prec x$ and computing the sought comprehensive Gr\"obner system using the algorithm {\sc PGB} give the results displayed in the following table:
\[
\footnotesize\aligned
 \begin{tabular}{|c||c||c|}
 \cline{1-3}
$E_i$&$N_i$&$G_i$\\
 \cline{1-3}
$\{a,b,c\}$ & $\{\  \}$ & $\{x,y\}$ \\
\cline{1-3}
$\{a, b \}$&$\{c\}$ & $\{cx^2-y, cy^2-x\}$\\
\cline{1-3}
$\{a^6-b^6, a^3c-b^3, b^3c-a^3,$&$\{b\}$&$ \{bx-acy, by-a\}$\\
 \ \ \ $ac^2-a, bc^2-b\}$&&\\
 \cline{1-3}
$\{\  \}$&$\{a^6-b^6, a^3c-b^3, b^3c-a^3,$ & $\{1\}$\\
 & $ac^2-a, bc^2-b\}$&\\
\cline{1-3}
 \end{tabular}
 \endaligned
\]
Accordingly, the algorithm divides the solution set of the system $\Sigma$ into four partitions, each of them corresponds to one of the above rows. Let us explain what each of these rows means. For the first row, we have the null conditions $E_1=\{a,b,c\}$ and there is no any non-null condition. This means that if the elements of $E_1$ are null, namely if $a=b=c=0$ then the system $\Sigma$ reduces to the system $G_1=\{x=0,y=0\}$ which obviously has the single solution $(0,0)$. For the second row, we have null conditions $E_2=\{a,b\}$ and non-null condition $N_2=\{c\}$. This means that if $a=b=0$ and $c\neq 0$ then the system $\Sigma$ reduces to $G_2=\{cx^2-y=0,cy^2-x=0\}$ which has the solution set $\{(\frac{1}{c},\frac{1}{c}), c\in\mathbb C\}\cup\{(0,0),(\frac{-1+\sqrt{3}i}{2c},\frac{-1-\sqrt{3}i}{2c}),(\frac{-1-\sqrt{3}i}{2c},\frac{-1+\sqrt{3}i}{2c})\}$. Similar interpretation holds for the third row. Finally, the last row means that if none of the previous null conditions holds, namely if $E_1,E_2,E_3\neq 0$, then the under consideration system $\Sigma$ reduces to $G_4=\{1=0\}$ which of course its solution set is empty.
\end{Example}

\subsection{Implementation}

In \cite{Implementation}, we have implemented the designed algorithm in {\sc Maple} 15 which is available online as a library entitled {\sc CRAut} (together with a sample file relevant to the following example). To do it, at first we implemented the recent algorithm {\sc PGB} of Kapur, Sun and Wang \cite{kapur2}. By means of this auxiliary algorithm, we employ techniques of comprehensive Gr\"obner systems to consider and solve the appearing parametric {\it linear} systems {\sf Sys} in the parametric case. Our implementation {\sc CRAut} enables one to compute the desired algebras of infinitesimal CR-automorphisms associated to homogeneous and weighted homogeneous CR-manifolds in two non-parametric and parametric cases and also to evaluate the effectiveness of the algorithm designed in this paper.

\begin{Example}
({\it reverification of the Beloshapka's conjecture in the lengths $\rho=1,\ldots,5$, cf.} \cite{Beloshapka2012, Mamai2009, Shananina2000}).
Consider the following weighted homogeneous rigid defining equations:
\begin{equation}
\label{Belo-Models}
\footnotesize
\aligned
w_1-\overline w_1&=2i\,z\overline z,
\\
 w_2-\overline w_2&=2i\,\big(z^2\overline{z}
+
\overline{z}^2z\big), \ \ \ \  w_3-\overline w_3=2\,\big(z^2\overline{z}
-\overline{z}^2\,z\big),
\\
w_4-\overline w_4&=2i\,\big(z^3\overline{z}
+
\overline{z}^3z\big)+2i{\bf a}\,z^2\overline{z}^2, \ \ \ \ w_5-\overline w_5=2\,\big(z^3\overline{z}
-\overline{z}^3z\big)+2i{\bf b}\,z^2\overline z^2, \ \ \ \ w_6-\overline w_6=2i\,z^2\overline{z}^2,  \ \ \ \ \ \ \ \ a\in\mathbb R,
\\
w_7-\overline w_7&=2i\,\big(z^4\overline z+z\overline z^4\big)+{\bf c}i\,\big(w_1+\overline w_1\big)\,\big(z^2\overline z+z\overline z^2\big),
 \ \ \ \ w_8-\overline w_8=2\big(z^4\overline z-z\overline z^4\big)+{\bf d}i\,\big(w_1+\overline w_1\big)\,\big(z^2\overline z+z\overline z^2\big),
 \\
 w_9-\overline w_9&=2i\,\big(z^3\overline z^2+z^2\overline z^3\big)+{\bf e}i\,\big(w_1+\overline w_1\big)\,\big(z^2\overline z+z\overline z^2\big),
\\
w_{10}-\overline w_{10}&=2\,\big(z^3\overline z^2-z^2\overline z^3\big)+{\bf f}i\,\big(w_1+\overline w_1\big)\,\big(z^2\overline z+z\overline z^2\big),
\\
w_{11}-\overline w_{11}&=\big(w_1+\overline w_1\big)\,\big(z^2\overline z-z\overline z^2\big), \ \ \ \  w_{12}-\overline w_{12}=\big(w_1+\overline w_1\big)\,\big(z^2\overline z+z\overline z^2\big), \ \ \ \ \ \ \ \ {\bf a,b,c,d,e,f\in \mathbb R},
\endaligned
\end{equation}
and let $M_{\sf k}$ be the Beloshapka's CR-model of CR-dimension 1 and codimension ${\sf k}=1,\ldots, 12$, represented in coordinates $(z,w_1,\ldots,w_{\sf k})$ in $\mathbb C^{{\sf k}+1}$.  For ${\sf k}=1,\ldots,5$, $M_{\sf k}$ is represented as the graph of the above first $\sf k$ equations. For ${\sf k}=6,\ldots,11$, it is represented again by the first $\sf k$ equations but with the assumption ${\bf a,b}=0$. Finally, $M_{\sf 12}$ is represented as the graph of the above 12 equations with the assumption that all the appearing six parameters $\bf a,b,c,d,e,f$ are vanished\,\,---\,\,namely a non-parametric model as $M_{\sf 1},M_{\sf 2},M_{\sf 3}$ and $M_{\sf 6}$ are.  These twelve CR-manifolds encompass all the Beloshapka's models up to the length five and are constructed by Shananina and Mamai \cite{Mamai2009, Shananina2000}. They also computed the associated Lie algebras of infinitesimal CR-automorphisms, thought Mamai did not present the outputs, perhaps because of the length of them. It is also known that these models are all homogeneous (\cite{Beloshapka2004}). By means of our implementation, we have computed the associated Lie algebras of infinitesimal CR-automorphisms. The following table displays some properties of the obtained results, where the timings were conducted on a personal laptop with Intel(R) Core(TM) {\it i}7 CPU@2.80 GHz and 6.00 GB of RAM:

\begin{equation*}
\aligned
\begin{tabular}{|c||c|c|c|c|c|c|c|c|c|c|c|c|}
  \hline
  {\rm Model} & $M_{\sf 1}$ & $M_{\sf 2}$ & $M_{\sf 3}$ & $M_{\sf 4}$ & $M_{\sf 5}$ & $M_{\sf 6}$ & $M_{\sf 7}$ & $M_{\sf 8}$ & $M_{\sf 9}$ & $M_{\sf 10}$ & $M_{\sf 11}$ & $M_{\sf 12}$\\
  \hline
  {\rm Time (sec.)} & 0.5 & 0.5 & 1 & 3 & 6.2 & 5.8 & 17.6 & 52.2 & 131.5 & 340 & 198 & 22 \\
  \hline
  $\rho$ & 2 & 3 & 3 & 4 & 4 & 4 & 5 & 5 & 5 & 5 & 5 & 5 \\
  \hline
  $\varrho$ & 2 & 0 & 0 & 0 & 0 & 0 & 0 & 0 & 0 & 0 & -1 & -1\\
  \hline
  {\rm dim.} & 8 & 5 & 7 & 7 & 9\, \vline\, 8 & 10 & 10\, \vline\, 9 & 12\, \vline\, 10 & 12\, \vline\, 11 & 14\, \vline\, 12 & 13 & 14 \\
  \hline
\end{tabular}
\endaligned
\end{equation*}
The last row of the above table needs some explanation. In fact for some models, the dimension of the associated algebra is not unique and depends on the values of the extant parameters. Being more precise, in this table we observe two different values for the dimensions associated to some models. In such cases, the left number is the dimension of the desired Lie algebra associated to the model whereas all the appearing parameters vanish identically; otherwise the dimension is equal to the number at the right hand side. For example for $M_{\sf 8}$, we have ${\rm dim}\big(\frak{aut}_{CR}(M_{\sf 8})\big)=12$ if ${\bf c,d}=0$ and $=10$, otherwise. More precisely, in the case that ${\bf c,d}\neq 0$, we have the basis elements of the components of $\frak{aut}_{CR}(M_{\sf 8})$ as:
\begin{equation}
\aligned
\label{k=8}
\frak g_{-5}&=\langle\partial_{w_7},\partial_{w_8}\rangle, \ \ \ \ \ \frak g_{-4}=\langle\partial_{w_6},\partial_{w_5},\partial_{w_4}\rangle, \ \ \ \ \ \frak g_{-3}=\langle\partial_{w_3},\partial_{w_2}\rangle,
\\
\frak g_{-2}&=\langle\partial_{w_1}+{\bf a}w_2\partial_{w_7}+{\bf b}w_2\partial_{w_8}\rangle,
\\
\frak g_{-1}&=\big\langle-2w_5\partial_{w_7}+2w_4\,\partial_{w_8}+z\partial_{w_1}+w_3\partial_{w_6}+z^3\partial_{w_4}+2w_1\partial_{w_3}+
z^2\partial_{w_2}+z^4\partial_{w_7}-\frac{3}{2}w_3\partial_{w_4}+\frac{3}{2}w_2\partial_{w_5}+
\\
&+{\bf c}z^2w_1\partial_{w_7}+
{\bf c}w_6\partial_{w_7}+{\bf d}z^2w_1\partial_{w_8}+{\bf d}w_6\partial_{w_8}+\frac{i}{2}\big(\partial_{z}-2z^2\partial_{w_3}-2z^3\partial_{w_5}-2z^4\partial_{w_8}\big),
\\
&{\bf c}z^2w_1\partial_{w_7}+{\bf d}z^3\partial_{w_4}+{\bf d}z^2w_1\partial_{w_8}+{\bf d}z\partial_{w_1}+{\bf d}z^4\partial_{w_7}+
{\bf d}z^2\partial_{w_2}-{\bf cb}w_1^2\partial_{w_7}-\frac{3}{2}{\bf d}w_3\partial_{w_5}-\frac{3}{2}{\bf d}w_2\partial_{w_4}
\\
&-
{\bf d}w_2\partial_{w_6}-{\bf d}z^4\partial_{w_8}-{\bf d}z^3\partial_{w_5}-{\bf d}z^2\partial_{w_3}-\frac{1}{2}{\bf d}\partial_{z}-
2{\bf d}w_4\partial_{w_7}-2{\bf d}w_1\partial_{w_2}-2{\bf d}w_5\partial_{w_8}-{\bf d}w_1^2\partial_{w_8}
\big\rangle,
\\
\frak g_0&=\langle\,\rangle.
\endaligned
\end{equation}
In the case that ${\bf c,d}=0$, the basis elements of $\frak g_i, i=-1,\ldots,-5$ are as above with of course ${\bf c,d}=0$ while in this case $\frak g_0$ has two basis elements as follows:
\[\aligned
\frak g_0&=\big\langle-w_8\partial_{w_7}+\frac{1}{3}w_2\partial_{w_3}+w_7\partial_{w_8}-\frac{1}{3}w_3\partial_{w_2}+\frac{2}{3}w_4\partial_{w_5}-
\frac{2}{3}w_5\partial_{w_4}+\frac{1}{3}iz\partial_{z},
\\
&\ \ \ \ \ \ \ \ \frac{1}{5}z\partial_{z}+\frac{3}{5}w_3\partial_{w_3}+\frac{3}{5}w_2\partial_{w_2}+\frac{4}{5}w_5\partial_{w_5}+
\frac{4}{5}w_6\partial_{w_6}+\frac{4}{5}w_4\partial_{w_4}+w_7\partial_{w_7}+w_8\partial_{w_8}+\frac{2}{5}w_1\partial_{w_1}\big\rangle.
\endaligned
\]
One easily concludes from the results of the above table,  all the above twelve models have rigidity, as is the main result of Mamai in \cite{Mamai2009}. Together with the library {\sc CRAut}, we also have put a sample file concerning the computations of this example.
\end{Example}

According to the above table, it took just 13 minutes from the implementation to verify the Beloshapka's conjecture in the lengths $\rho=1,\ldots,5$. This shows the effectiveness of the algorithm, designed in this paper.

\begin{Remark}
A glance on the above timings shows how it increases mostly the complexity of computations as one
passes from each model to the next by adding just one variable and one defining equation to the previous
ones. The following diagram may be helpful to compare the appearing timings. Moreover, actually the computations in the case of parametric defining equations are more complicated in comparison to those of the non-parametric case. For example, compare the timings corresponding to the models $M_{11}$ and $M_{\sf 12}$\,\,---\,\,notice that the defining equations of $M_{\sf 12}$ are non-parametric.
\end{Remark}

\begin{center}
\includegraphics[scale=0.65]{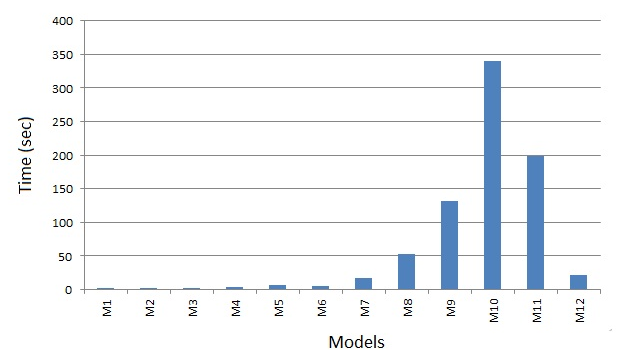}
\end{center}

\begin{Example}
Add the following rigid defining equations to the list \thetag{\ref{Belo-Models}}:
\begin{equation}
\label{Belo-Models-2}
\aligned
w_{13}-\overline w_{13}&=2i\,\big(z^5\overline{z}
+
\overline{z}^5z\big), \ \ \ \ w_{14}-\overline w_{14}=2\,\big(z^5\overline{z}
-
\overline{z}^4z\big)
, \ \ \ \ w_{15}-\overline w_{15}=2i\,\big(z^4\overline{z}^2
+
\overline{z}^4z^2\big),
\\
w_{16}-\overline w_{16}&=2\,\big(z^4\overline{z}^2
-\overline{z}^4z^2\big), \ \ \ \ w_{17}-\overline w_{17}=2i\,z^3\overline{z}^3,
\endaligned
\end{equation}
and for ${\sf k}=13,\ldots,17$ let $M_{\sf k}$ to be the CR-manifold of CR-dimension one and codimension $\sf k$ represented as the graph of the equations of \thetag{\ref{Belo-Models}} together the first ${\sf k}-12$ equations of the above list. These are the next five rigid Beloshapka's models which are the very first models of the length six. Here, we also compute\,\,---\,\,for the first time\,\,---\,\,the desired Lie algebras associated to these models and the results are displayed in the following table:
\begin{equation*}
\aligned
\begin{tabular}{|c||c|c|c|c|c|}
  \hline
  {\rm Model} & $M_{\sf 13}$ & $M_{\sf 14}$ & $M_{\sf 15}$ & $M_{\sf 16}$ & $M_{\sf 17}$ \\
  \hline
  {\rm Time (sec.)} & 83.5 & 152 & 286& 545  & 1157 \\
  \hline
  $\rho$ & 6 & 6 & 6 & 6 & 6 \\
  \hline
  $\varrho$ & -1 & -1 & -1 & -1 & -1 \\
  \hline
  {\rm dim} & 15 & 16 & 17 & 18 & 19 \\
  \hline
\end{tabular}
\endaligned
\end{equation*}
\end{Example}

\subsection*{Acknowledgment}

The authors gratefully acknowledge the helpful discussions of Valerii Beloshapka, Mauro Nacinovich and Andrea Spiro during the preparation of this paper. The research of the second author was in part supported by a grant from IPM (No. 92550420).

\end{document}